\documentclass[aap,MSNbibl,citesort,dvips]{arximspdf}
\usepackage{mathbh}
\usepackage{graphicx}

\doi{10.1214/10-AAP695}
\volume{21}
\issue{1}
\pubyear{2011}
\firstpage{283}
\lastpage{311}

\makeatletter
\newcommand{\suptwo}[2]{\mathop{\mathop{\sup}_{#1}}_{#2}} 
\newcommand{\inftwo}[2]{\mathop{\mathop{\inf}_{#1}}_{#2}} 
\newcommand{\cA}{\mathcal{A}}
\newcommand{\cF}{\mathcal{F}}
\newcommand{\cG}{\mathcal{G}}
\newcommand{\cH}{\mathcal{H}}
\newcommand{\cW}{\mathcal{W}}
\newcommand{\cX}{\mathcal{X}}

\newcommand{\io}{{\iota}}
\newcommand{\x}{{\bar Y}}
\newcommand{\xis}{\bar Y_{\iota(s)}}
\newcommand{\xit}{\bar Y_{\iota(t)}}
\newcommand{\xism}{\bar Y_{\iota(s-)}}
\newcommand{\xitm}{\bar Y_{\iota(t-)}}
\newcommand{\xes}{\bar Y_{\eta(s)}}
\newcommand{\xet}{\bar Y_{\eta(t)}}
\newcommand{\xesm}{\bar Y_{\eta(s-)}}
\newcommand{\xetm}{\bar Y_{\eta(t-)}}
\newcommand{\z}{{\bar\ups}}
\newcommand{\zis}{\bar\ups_{\iota(s)}}
\newcommand{\zit}{\bar\ups_{\iota(t)}}
\newcommand{\zism}{\bar\ups_{\iota(s-)}}
\newcommand{\zitm}{\bar\ups_{\iota(t-)}}
\newcommand{\zet}{\bar\ups_{\eta(t)}}
\newcommand{\zesm}{\bar\ups_{\eta(s-)}}
\newcommand{\zetm}{\bar\ups_{\eta(t-)}}
\newcommand{\ups}{\Upsilon}

\newcommand{\sfrac}[2]{\frac{#1}{#2}}

\newcommand{\ind}{\mathbh{1}}
\newcommand{\dd}{\mathrm{d}}
\renewcommand{\th}{\theta}
\newcommand{\eps}{\varepsilon}
\newcommand{\Sig}{\Sigma}
\newcommand{\vphi}{\varphi}
\newcommand{\Lam}{\Lambda}
\newcommand{\lfl}{\lfloor}
\newcommand{\rfl}{\rfloor}

\newcommand{\cov}{\operatorname{cov}}
\newcommand{\var}{\operatorname{var}}
\newcommand{\IP}{\mathbb{P}}
\newcommand{\II}{\mathbb{I}}
\newcommand{\IIJ}{\mathbb{J}}
\newcommand{\IN}{\mathbb{N}}
\newcommand{\IZ}{\mathbb{Z}}
\newcommand{\IC}{\mathbb{C}}
\newcommand{\IR}{\mathbb{R}}
\newcommand{\IE}{\mathbb{E}}
\newcommand{\vth}{\vartheta}

\newtheorem{theorem}{Theorem}[section]
\newtheorem{prop}[theorem]{Proposition}
\newtheorem{lemma}[theorem]{Lemma}
\newtheorem{cor}[theorem]{Corollary}
\newproclaim{assu}{Assumption}
\DeclareMathAlphabet\mathcaligr{OMS}{cmsy}{m}{n}
\renewcommand{\mathcal}{\mathcaligr}
\makeatother

\begin{document}
\begin{frontmatter}

\title{Multilevel Monte Carlo algorithms for L\'evy-driven SDEs with
Gaussian correction}
\runtitle{MLMC for L\'evy SDE with Gaussian correction}

\begin{aug}
\author{\fnms{Steffen} \snm{Dereich}\corref{}\ead[label=e1]{dereich@mathematik.uni-marburg.de}
\ead[label=u1,url]{http://www.mathematik.uni-marburg.de/\texttildelow dereich}}
\runauthor{S. Dereich}
\affiliation{Philipps-Universit\"at Marburg}
\address{Philipps-Universit\"at Marburg\\
Fb. 12, Mathematik und Informatik\\
Hans-Meerwein-Stra\ss e\\
D-35032 Marburg\\
Germany\\
\printead{e1}\\
\printead{u1}}
\end{aug}

\received{\smonth{8} \syear{2009}}
\revised{\smonth{1} \syear{2010}}

%
\begin{abstract}
We introduce and analyze multilevel Monte Carlo algorithms for the
computation of $\IE f(Y)$, where
$Y=(Y_t)_{t\in[0,1]}$ is the solution of a multidimensional L\'
evy-driven stochastic differential
equation and $f$ is a real-valued function on the path space. The
algorithm relies on approximations
obtained by simulating large jumps of the L\'evy process individually
and applying a Gaussian approximation
for the small jump part.
Upper bounds are provided for the worst case error over the class of
all measurable real functions $f$
that are Lipschitz continuous with respect to the supremum norm. These
upper bounds are easily tractable
once one knows the behavior of the L\'evy measure around zero.

In particular, one can derive upper bounds from the Blumenthal--Getoor
index of the L\'evy process.
In the case where the Blumenthal--Getoor index is larger than one, this
approach is superior to
algorithms that do not apply a Gaussian approximation.
If the L\'evy process does not incorporate a Wiener process or if the
Blumenthal--Getoor index
$\beta$ is larger than $\frac43$, then the upper bound is of order
$\tau^{- ({4-\beta})/({6\beta})}$
when the runtime $\tau$ tends to infinity. Whereas in the case, where
$\beta$ is in $[1,\frac43]$ and
the L\'evy process has a Gaussian component, we obtain bounds of order
$\tau^{- \beta/(6\beta-4)}$.
In particular, the error is at most of order $\tau^{- 1/6}$.
\end{abstract}

%

\end{frontmatter}

\section{Introduction}\label{intro}

Let $d_Y\in\IN$ and denote by $D[0,1]$ the Skorokhod space of
functions mapping $[0,1]$ to $\IR^{d_Y}$ endowed with its
Borel-$\sigma$-field.
In this article, we analyze numerical schemes for the evaluation of
\[
S(f):=\IE[ f(Y)],
\]
where
\begin{itemize}
\item
$Y=(Y_t)_{t\in[0,1]}$ is a solution to a multivariate stochastic
differential equation driven by a multidimensional L\'evy process (with
state space $\IR^{d_Y}$), and
\item$f\dvtx D[0,1]\to\IR$ is a Borel measurable function that is
Lipschitz continuous with respect to the \textit{supremum} norm.
\end{itemize}

This is a classical problem which appears for instance in finance,
where $Y$ models the risk neutral stock price and $f$ denotes the
payoff of a (possibly path dependent) option, and in the past several
concepts have been employed for dealing with it.

A common stochastic approach is to perform a Monte Carlo simulation of
numerical approximations to the solution $Y$.
Typically, the Euler or Milstein schemes are used to obtain
approximations. Also higher order schemes can be applied provided that
samples of iterated It\^o integrals are supplied and the coefficients
of the equation are sufficiently regular. In general, the problem is
tightly related to \textit{weak approximation} which is, for instance,
extensively studied in the monograph by Kloeden and Platen \cite
{KloePla92} for diffusions.

Essentially, one distinguishes between two cases. Either $f(Y)$ depends
only on the state of $Y$ at a fixed time or alternatively it depends on
the whole trajectory of~$Y$. In the former case, extrapolation
techniques can often be applied to increase the order of convergence,
see \cite{TalTub90}.
For L\'evy-driven stochastic differential equations, the Euler scheme
was analyzed in \cite{ProTal97} under the assumption that the
increments of the L\'evy process are simulatable. Approximate
simulations of the L\'evy increments are considered in \cite{JKMP05}.

In this article, we consider functionals $f$ that depend on the whole
trajectory. 
Concerning results for diffusions, we refer the reader to the
monograph \cite{KloePla92}. For L\'evy-driven stochastic differential
equations, limit theorems in distribution are provided in \cite{Jac04}
and \cite{Rub03} for the discrepancy between the genuine solution and
Euler approximations.

Recently, Giles \cite{Gil08,Gil08b} (see also \cite{Hei98})
introduced the so-called \textit{multilevel Monte Carlo method} to
compute $S(f)$. It is very efficient when $Y$ is a diffusion. Indeed,
it even can be shown that it is---in some sense---optimal, see \cite
{CDMR08}. For L\'evy-driven stochastic differential equations,
multilevel Monte Carlo algorithms are first introduced and studied
in \cite{DerHei10}. Let us explain their findings in terms of the
Blumenthal--Getoor index (BG-index) of the driving L\'evy process which
is an index in $[0,2]$. It measures the frequency of small jumps,
see (\ref{eq0826-1}), where a large index corresponds to a process
which has small jumps at high frequencies. In particular, all L\'evy
processes which have a finite number of jumps has BG-index zero.
Whenever the BG-index is smaller or equal to one, the algorithms of
\cite{DerHei10} have worst case errors at most of order $\tau^{-
1/2}$, when the runtime $\tau$ tends to infinity. Unfortunately, the
efficiency decreases significantly for larger Blumenthal--Getoor indices.

Typically, it is not feasible to simulate the increments of the L\'evy
process perfectly, and one needs to work with approximations. This
necessity typically worsens the performance of an algorithm, when the
BG-index is larger than one due to the higher frequency of small jumps.
It represents the main bottleneck in the simulation.
In this article, we consider approximative L\'evy increments that
simulate the large jumps and approximate the small ones by a normal
distribution (\textit{Gaussian approximation}) in the spirit of Asmussen
and Rosi{\'n}ski \cite{AsmRos01} (see also \cite{CohRos07}).
Whenever the BG-index is larger than one, this approach is superior to
the approach taken in \cite{DerHei10}, which neglects small jumps in
the simulation of L\'evy increments.

\begin{figure}

\includegraphics{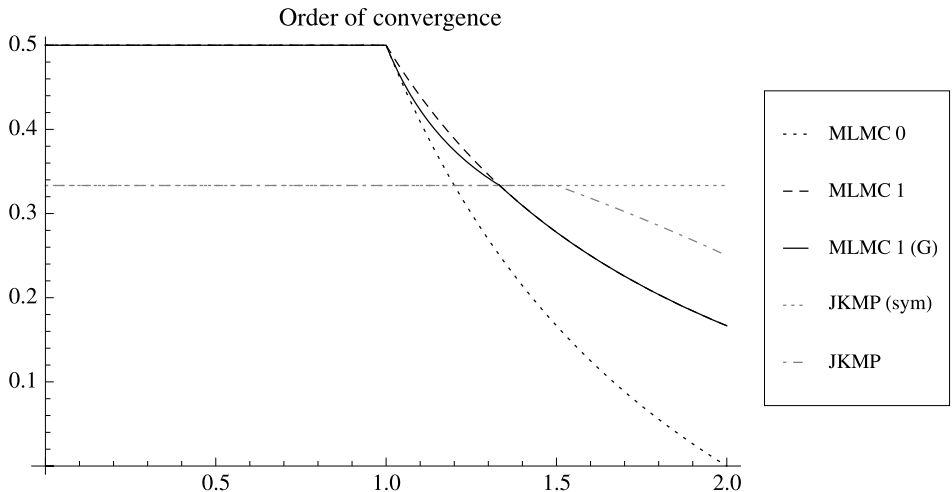}

\caption{Order of convergence in dependence on the Blumenthal--Getoor index.}\label{fig1}
\end{figure}

To be more precise, we establish a new estimate for the Wasserstein
metric between an approximative solution with Gaussian approximation
and the genuine solution, see Theorem \ref{th0310-1}. It is based on a
consequence of Zaitsev's generalization \cite{Zai98} of the Koml\'os--Major--Tusn\'ady coupling
\cite{KMT75,KMT76} which might be of its
own interest itself, see Theorem~\ref{th0217-1}.
With these new estimates, we analyze a class of multilevel Monte Carlo
algorithms together with a \textit{cost} function which measures the
computational complexity of the individual algorithms.
We provide upper error bounds for individual algorithms and optimize
the error over the parameters under a given cost constraint. When the
BG-index is larger than one, appropriately adjusted algorithms lead to
significantly smaller \textit{worst case errors} over the class of
Lipschitz functionals than the ones analyzed so far, see Theorem~\ref
{th:le}, Corollary~\ref{cor1}  and Figure~\ref{fig1}. In particular,
one always obtains numerical schemes with errors at most of order $\tau
^{-1/6}$ when the runtime $\tau$ of the algorithm tends to infinity.


\subsection*{Notation and universal assumptions}

We denote by $|\cdot|$ the Euclidean norm for vectors as well as the
Frobenius norm for matrices and let $\|\cdot\|$ denote the supremum
norm over the interval $[0,1]$. $X=(X_t)_{t\ge0}$ denotes an
$d_X$-dimensional $L^2$-integrable L\'evy process. By the L\'evy--Khintchine formula, it is characterized by a square integrable L\'
evy-measure $\nu$ [a Borel measure on $\IR^{d_X}\backslash\{0\}$
with $\int|x|^2 \nu(\dd x)<\infty$], a positive semi-definite
matrix $\Sig\Sig^*$ ($\Sig$ being a $d_X\times d_X$-matrix), and a
drift $b\in\IR^{d_X}$ via
\[
\IE e^{i \langle\th, X_t\rangle}=e^{t \psi(\th)},
\]
where
\[
\psi(\th)= \frac12 |\Sig^* \th|^2 + \langle b, \th\rangle+ \int
_{\IR^{d_X}}  \bigl(e^{i\langle\th, x\rangle}-1- i\langle\th
,x\rangle \bigr) \nu(\dd x).
\]
Briefly, we call $X$ a $(\nu, \Sig\Sig^*, b)$-L\'evy process, and
when $b=0$, a $(\nu,\Sig\Sig^*)$-L\'evy martingale. All L\'evy
processes under consideration are assumed to be c\`adl\`ag.
As is well known, we can represent $X$ as sum of three independent processes
\[
X_t= \Sig W_t + L_t+ b t,
\]
where $W=(W_t)_{t\ge0}$ is a $d_X$-dimensional Wiener process and
$L=(L_t)_{t\ge0}$ is a $L^2$-martingale that comprises the compensated
jumps of $X$.
We consider the integral equation
%
\begin{equation}\label{SDE}
Y_t =y_0+\int_0^t a(Y_{t-}) \,\dd X_t,
\end{equation}
where $y_0\in\IR^{d_Y}$ is a fixed deterministic initial value.
We impose the standard Lipschitz assumption on the function $a\dvtx \IR
^{d_Y}\to\IR^{d_Y\times d_X}$: for a fixed $K<\infty$, and all
$y,y'\in\IR^{d_Y}$, one has
\[
|a(y)-a(y')|\leq K |y-y'| \quad  \mbox{and}\quad   |a(y_0)|\leq K.
\]
Furthermore, we assume without further mentioning that
\[
\int|x|^2 \nu(\dd x)\leq K^2 , \qquad  |\Sig|\le K \quad \mbox{and}\quad
|b|\le K .
\]
We refer to the monographs \cite{Ber98} and \cite{Sat99} for details
concerning L\'evy processes. Moreover, a comprehensive introduction to
the stochastic calculus for discontinuous semimartingales and, in
particular, L\'evy processes can be found in \cite{Pro05} and~\cite{Apple04}.

In order to approximate the small jumps of the L\'evy process, we need
to impose a uniform ellipticity assumption.

\renewcommand{\theassu}{UE}
\begin{assu} \label{assu1}
There are $\mathfrak{h}\in(0,1]$, $\vth\ge1$ and a linear subspace
$\cH$ of $\IR^{d_X}$ such that for all $h\in(0,\mathfrak{h}]$ the
L\'evy measure $\nu|_{B(0,h)}$ is supported on $\cH$ and satisfies
\[
\frac1\vth\int_{B(0,h)} \langle y,x\rangle^2  \nu(\dd x) \le\int
_{B(0,h)} \langle y',x\rangle^2  \nu(\dd x) \le\vth\int_{B(0,h)}
\langle y,x\rangle^2  \nu(\dd x)
\]
for all $y,y'\in\cH$ with $|y|=|y'|$.
\end{assu}

\subsection*{Main results}

We consider a class of multilevel Monte Carlo algorithms $\cA$
together with a cost function $\operatorname{cost} \dvtx \cA\to[0,\infty)$
that are introduced explicitly in Section~\ref{sec:alg}. For each
algorithm $\widehat S\in\cA$, we denote by $\widehat S(f)$ a
real-valued random variable representing the random output of the
algorithm when applied to a given measurable function $f\dvtx D[0,1]\to\IR
$. We work in the real number model of computation, which means that we
assume that arithmetic operations with real numbers and comparisons can
be done in one time unit, see also \cite{Nov95}. 
Our cost function represents the runtime of the algorithm reasonably
well when supposing that
\begin{itemize}
\item one can sample from the distribution $\nu|_{B(0,h)^c}/\nu
(B(0,h)^c)$ and the uniform distribution on $[0,1]$ in constant time,
\item one can evaluate $a$ at any point $y\in\IR^{d_Y}$ in constant
time, and
\item$f$ can be evaluated for piecewise constant functions in less
than a constant multiple of its breakpoints plus one time units.
\end{itemize}
As pointed out below, in that case, the average runtime to evaluate
$\widehat S(f)$ is less than a constant multiple of $\operatorname
{cost}(\widehat S)$. We analyze the minimal \textit{worst case error}
\[
\operatorname{err}(\tau)= \inftwo{\widehat S\in\cA:}{\operatorname
{cost}(\widehat S)\le\tau} \sup_{f\in\operatorname{Lip}(1)} \IE
[|S(f)-\widehat S(f)|^2]^{1/2},\qquad \tau\ge1.
\]
Here and elsewhere, $\operatorname{Lip}(1)$ denotes the class of measurable
functions $f\dvtx D[0,\break 1]\to\IR$ that are Lipschitz continuous with respect
to \textit{supremum} norm with coefficient one.

In this article, we use asymptotic comparisons. We write $f\approx g$ for
$ 0< \liminf\frac fg \le\limsup\frac fg <\infty, $
and $f\precsim g$ or, equivalently $g\succsim f$, for
$ \limsup\frac fg <\infty. $
Our main findings are summarized in the following theorem.

\begin{theorem}\label{th:le}
Assume that Assumption  \ref{assu1}
is valid and
let $g\dvtx (0,\infty)\to(0,\infty)$ be a decreasing and invertible
function such that for all $h>0$
\[
\int\frac{|x|^2}{h^2}\wedge1  \nu(\dd x)\le g(h)
\]
and, for a fixed $\gamma>1$,
%
\begin{equation}\label{eq1126-1}
g \biggl(\frac\gamma2 h \biggr)\ge2 g(h)
\end{equation}
for all sufficiently small $h>0$.
\begin{enumerate}[(II)]
\item[(I)]
If $\Sig=0$ or
\[
g^{-1}(x)\succsim x^{-3/4}  \qquad\mbox{as }  x\to\infty,
\]
then
\[
\operatorname{err}(\tau)\precsim g^{-1} ((\tau\log\tau)^{2/3}
) \tau^{1/6} (\log\tau)^{2/3}  \qquad\mbox{as }  \tau\to\infty.
\]
\item[(II)] If
\[
g^{-1}(x)\precsim x^{-3/4}  \qquad\mbox{as }  x\to\infty,
\]
then
\[
\operatorname{err}(\tau)\precsim\sqrt{\frac{\log\tau}{g^*(\tau)}}
\qquad\mbox{as }  \tau\to\infty,
\]
where $g^*(\tau)=\inf\{x>1\dvtx  x^3 g^{-1}(x)^2 (\log x)^{-1}\ge\tau\}$.
\end{enumerate}
\end{theorem}

The class of algorithms $\cA$ together with appropriate parameters
which establish the error estimates above are stated explicitly in
Section~\ref{sec:alg}.

In terms of the Blumenthal--Getoor index
%
\begin{equation}\label{eq0826-1}
\beta:=\inf \biggl\{p>0\dvtx  \int_{B(0,1)} |x|^p  \nu(\dd x) < \infty
 \biggr\}\in[0,2]
\end{equation}
we get the following corollary.

\begin{cor} \label{cor1} Assume that Assumption  \ref{assu1}
is valid and that the BG-index satisfies $\beta\ge1$.
If $\Sig=0$ or $\beta\ge\frac43$, then
\[
\sup\{\gamma\ge0\dvtx  \operatorname{err}(\tau)\precsim\tau^{-\gamma}\}\ge
\frac{4-\beta}{6\beta},
\]
and, if $\Sig\not=0$ and $\beta< \frac43$,
\[
\sup\{\gamma\ge0\dvtx  \operatorname{err}(\tau)\precsim\tau^{-\gamma}\}\ge
\frac{\beta}{6\beta-4}.
\]
\end{cor}

\subsubsection*{Visualization of the results and relationship to other work}

Figure~\ref{fig1} illustrates our findings and related results. The
$x$-axis and $y$-axis represent the Blumenthal--Getoor index and the
order of convergence, respectively. Note that MLMC 0 stands for the
multilevel Monte Carlo algorithm which does not apply a Gaussian
approximation, see \cite{DerHei10}. Both lines marked as MLMC 1
illustrate Corollary \ref{cor1}, where the additional (G) refers to
the case where the SDE comprises a Wiener process.

These results are to be compared with the results of Jacod et al. \cite
{JKMP05}. Here an approximate Euler method is analyzed by means of weak
approximation. In contrast to our investigation, the object of that
article is to compute $\IE f(X_T)$ for a fixed time $T>0$. Under quite
strong assumptions (for instance, $a$ and $f$ have to be four times
continuously differentiable and the eights moment of the L\'evy process
needs to be finite), they provide error bounds for a numerical scheme
which is based on Monte Carlo simulation of one approximative solution.
In the figure, the two lines quoted as JKMP represent the order of
convergence for general, respectively pseudo symmetrical, L\'evy processes.
Additionally to the illustrated schemes, \cite{JKMP05} provide an
expansion which admits a Romberg extrapolation under additional assumptions.

We stress the fact that our analysis is applicable to general path
dependent functionals and that our error criterion is the worst case
error over the class of Lipschitz continuous functionals with respect
to supremum norm. In particular, our class contains most of the
continuous payoffs appearing in finance.

We remark that our results provide upper bounds for the inferred error
and so far no lower bounds are known. The worst exponent appearing in
our estimates is $\frac16$ which we obtain for L\'evy processes with
Blumenthal--Getoor index $2$. Interestingly, this is also the worst
exponent appearing in \cite{RubWik03} in the context of strong
approximation of SDEs driven by subordinated L\'evy processes.

\subsubsection*{Agenda}

The article is organized as follows. In Section~\ref{sec:alg}, we
introduce a class of multilevel Monte Carlo algorithms together with a
cost function. Here, we also provide the crucial estimate for the mean
squared error which motivates the consideration of the Wasserstein
distance between an approximative and the genuine solution, see~(\ref
{eq0311-1}). Section~\ref{sec:3} states the central estimate for the
former Wasserstein distance, see Theorem~\ref{th0310-1}. In this
section, we explain the strategy of the proof and the structure of the
remaining article in detail. For the proof, we couple the driving L\'
evy process with a L\'evy process constituted by the large jumps plus a
Gaussian compensation of the small jumps and we write the difference
between the approximative and the genuine solution as a telescoping sum
including further auxiliary processes, see (\ref{eq0218-3}) and (\ref
{eq0218-4}). The individual errors are then controlled in Sections~\ref
{sec:4} and \ref{sec:5} for the terms which do not depend on the
particular choice of the coupling and in Section~\ref{sec:gauss} for
the error terms that do depend on the particular choice. In between, in
Section~\ref{sec:KMT}, we establish the crucial KMT like coupling
result for the L\'evy process.
Finally, in Section~\ref{sec:proofs}, we combine the approximation
result for the Wasserstein metric (Theorem~\ref{th0310-1}) with
estimates for strong approximation of stochastic differential equations
from \cite{DerHei10} to prove the main results stated above.

\section{Multilevel Monte Carlo}\label{sec:alg}

Based on a number of parameters, we define a multilevel Monte Carlo
algorithm $\widehat S$: We denote by $m$ and $n_1,\ldots,n_{m}$ natural
numbers and let $\eps_1,\ldots,\eps_m$ and $h_1,\ldots,h_m$ denote
decreasing sequences of positive reals.
Formally, the algorithm $\widehat S$ can be represented as a tuple
constituted by these parameters, and we denote by $\cA$ the set of all
possible choices for $\widehat S$.
We continue with defining processes that depend on the latter
parameters. For ease of notation, the parameters are omitted in the
definitions below.

We choose a square matrix $\Sig^{({m})}$ such that\vspace*{-1pt} $(\Sig^{({m})}
(\Sig^{({m})})^*)_{i,j}=\int_{B(0,h_m)} x_ix_j\times\break  \nu(\dd x)$.
Moreover, for $k=1,\ldots, m$, we let $L^{({k})}=(L^{({k})}_t)_{t\ge
0}$ denote the $(\nu|_{B(0,h_k)^c},\break 0)$-L\'evy martingale
which comprises the compensated jumps of $L$ that are larger than
$h_k$, that is
\[
L_t^{({k})} =\sum_{s\le t} \ind_{\{|\Delta L_s|\ge h_k\}} \Delta
L_s -t \int_{B(0,h_k)^c} x \nu(\dd x).
\]
Here and elsewhere, we denote $\Delta L_t=L_t-L_{t-}$.
We let $B=(B_t)_{t\ge0}$ be an independent Wiener process (independent
of $W$ and $L^{({k})}$), and consider, for $k=1,\ldots,m$, the
processes $\cX^{({k})} =(\Sig W_t +\Sig^{({m})} B_t + L^{({k})}_t+bt
)_{t\ge0}$ as driving processes. Let $\ups^{({k})}$ denote
the solution to
\[
\ups^{({k})}_t= y_0+\int_0^t a\bigl(\ups^{({k})}_{s-}\bigr)\, \dd\cX
_{\io^{({k})}(s)},
\]
where $(\io^{({k})}(t))_{t\ge0}$ is given via $\io^{({k})}(t)=
\max(\II^{({k})}\cap[0,t])$ and the set $\II^{({k})}$ is
constituted by the random times $(T^{({k})}_j)_{j\in\IZ_+}$ that
are inductively defined via $T^{({k})}_0=0$ and
\[
T^{({k})}_{j+1}=\inf\bigl\{ t\in\bigl( T^{({k})}_j,\infty\bigr)\dvtx  |\Delta
L_t|\ge h_k\mbox{ or } t= T^{({k})}_j+\eps_k\bigr\}.
\]
Clearly, $\ups^{({k})}$ is constant on each interval
$[T^{({k})}_j,T^{({k})}_{j+1})$ and one has
%
\begin{equation}\label{eq0514-3}
\ups^{({k})}_{T^{(k)}_{j+1}}=\ups^{({k})}_{T^{(k)}_{j}}+a\bigl(\ups
^{({k})}_{T^{(k)}_{j}}\bigr)  \bigl(\cX_{T^{(k)}_{j+1}}-\cX_{T^{(k)}_{j}}\bigr).
\end{equation}

Note that we can write
\[
\IE\bigl[f\bigl(\ups^{({m})}\bigr)\bigr]= \sum_{k=2}^{m} \IE\bigl[f\bigl( \ups^{({k})}\bigr)-f\bigl(\ups
^{({k-1})}\bigr)\bigr] +\IE\bigl[f\bigl(\ups^{({1})}\bigr)\bigr].
\]
The multilevel Monte Carlo algorithm---identified with $\widehat S$---estimates each expectation $\IE[f( \ups^{({k})})-f(\ups^{({k-1})})]$
(resp., $\IE[f(\ups^{({1})})]$) individually by sampling
independently $n_{k}$ (resp., $n_1$) versions of $f( \ups
^{({k})})-f(\ups^{({k-1})})$ [$f(\ups^{({1})})$] and taking the
average. The output of the algorithm is then the sum of the individual
estimates. We denote by $\widehat S(f)$ a random variable that models
the random output of the algorithm when applied to $f$.

\subsection*{The mean squared error of an algorithm}

The Monte Carlo algorithm introduced above induces the mean squared error
%
\begin{eqnarray*}
\operatorname{mse}(\widehat S,f)&=& \bigl|\IE[f(Y)]- \IE\bigl[f\bigl(\ups^{({m})}\bigr)\bigr]\bigr| ^2
+ \sum_{k=2}^{m} \frac1{n_k} \var\bigl(f\bigl( \ups^{({k})}\bigr)-f\bigl(\ups
^{({k-1})}\bigr)\bigr)\\
&&{} +\frac1{n_1} \var\bigl(f\bigl( \ups^{({1})}\bigr)\bigr),
\end{eqnarray*}
when applied to $f$.
For two $D[0,1]$-valued random elements $Z^{({1})}$ and $Z^{({2})}$, we
denote by $\cW(Z^{({1})}, Z^{({2})})$ the Wasserstein
metric of second-order with respect to supremum norm, that is
%
\begin{equation}\label{wasser}
\cW\bigl(Z^{({1})}, Z^{({2})}\bigr)=\inf_{\xi}  \biggl(\int\bigl\|
z^{({1})}-z^{({2})}\bigr\|^2\,\dd\xi\bigl(z^{({1})},z^{({2})}\bigr) \biggr)^{1/2},
\end{equation}
where the infimum is taken over all probability measures $\xi$ on
$D[0,1]\times D[0,1]$ having first marginal $\IP_{Z^{(1)}}$ and second
marginal $\IP_{Z^{(2)}}$. Clearly, the Wasserstein distance depends
only on the distributions of $Z^{({1})}$ and $Z^{({2})}$.
Now, we get for $f\in\operatorname{Lip}(1)$, that
%
\begin{eqnarray}\label{eq0311-1}
\operatorname{mse}\bigl(\widehat S,f\bigr)&\le&\cW\bigl(Y,\ups^{({m})}\bigr)^2 + \sum
_{k=2}^{m} \frac1{n_k} \IE\bigl[\bigl\| \ups^{({k})}-\ups^{({k-1})}\bigr\|
^2\bigr]\nonumber
\\[-8pt]
\\[-8pt]
&&{} +\frac1{n_1} \IE\bigl[\bigl\| \ups^{({1})}-y_0\bigr\|^2\bigr].\nonumber
\end{eqnarray}
We set
\[
\operatorname{mse}(\widehat S)=\sup_{f\in\operatorname{Lip}(1)} \operatorname
{mse}(\widehat S,f),
\]
and remark that estimate (\ref{eq0311-1}) remains valid for the worst
case error $\operatorname{mse}(\widehat S)$.

The main task of this article is to provide good estimates for the
Wasserstein metric $\cW(Y,\ups^{({m})})$. The remaining terms on
the right-hand side of (\ref{eq0311-1}) are controlled with estimates
from \cite{DerHei10}.

\subsection*{The cost function}

In order to simulate one pair $(\ups^{({k-1})},\ups^{({k})})$, we
need to simulate all displacements of $L$ of size larger or
equal to $h_{k}$ on the time interval $[0,1]$. Moreover, we need the
increments of the Wiener process on the time skeleton $(\II
^{({k-1})}\cup\II^{({k})})\cap[0,1]$. Then we can construct our
approximation via (\ref{eq0514-3}). In the real number model of
computation (under the assumptions described in the \hyperref[intro]{Introduction}), this
can be performed with runtime less than a multiple of the number of
entries in $\II^{({k})}\cap[0,1]$, see \cite{DerHei10} for a
detailed description of an implementation of a similar scheme. Since
\[
\IE \bigl[\# \bigl(\II^{({k})}\cap[0,1]\bigr) \bigr]\le1+\frac1{\eps
_{k}} +\IE\biggl[ \sum_{t\in[0,1]} \ind_{\{|\Delta L_t|\ge h_{k}\}}\biggr] =
\nu(B(0,h_{k})^c)+\frac1{\eps_{k}}+1,
\]
we define, for $\widehat S\in\cA$,
\[
\operatorname{cost}(\widehat S)=\sum_{k=1}^{m} n_k  \biggl[\nu
(B(0,h_{k})^c)+\frac1{\eps_{k}}+1 \biggr].
\]
Then supposing that $\eps_1\le1$ and $\nu(B(0,h_{k})^c)\le\frac
1{\eps_{k}}$ for $k=1,\ldots,m$, yields that
%
\begin{equation}\label{eq0514-2}
\operatorname{cost}(\widehat S)\le3 \sum_{k=1}^{m} n_k  \frac1{\eps_{k}}.
\end{equation}

\subsection*{Algorithms achieving the error rates of Theorem~\protect\ref{th:le}}

Let us now quote the choice of parameters which establish the error
rates of Theorem~\ref{th:le}.
In general, one chooses $\eps_k=2^{-k}$ and $h_k=g^{-1}(2^k)$ for
$k\in\IZ_+$. Moreover, in case (I), for sufficiently large $\tau$,
one picks
\begin{eqnarray}
m= \lfl\log_2 C_1 (\tau\log\tau)^{2/3} \rfl \quad \mbox{and}\quad
  n_k= \biggl\lfl C_2 \tau^{1/3} (\log\tau)^{-2/3} \frac
{g^{-1}(2^k)}{g^{-1}(2^m)} \biggr\rfl\nonumber\\
\eqntext{\mbox{for }k=1,\ldots,m,}
\end{eqnarray}
where $C_1$ and $C_2$ are appropriate constants that do not depend on
$\tau$.
In case (II), one chooses
\begin{eqnarray}
m= \lfl\log_2 C_1 g^*(\tau)  \rfl \quad \mbox{and}\quad
n_k= \biggl\lfl C_2 \frac{g^*(\tau)^2}{\log g^*(\tau)} \frac
{g^{-1}(2^k)}{g^{-1}(2^m)}\biggr \rfl \nonumber\\
 \eqntext{\mbox{ for }k=1,\ldots,m,}
\end{eqnarray}
where again $C_1$ and $C_2$ are appropriate constants. We refer the
reader to the proof of Theorem~\ref{th:le} for the error estimates of
this choice.

\section{Weak approximation}\label{sec:3}

In this section, we provide the central estimate for the Wasserstein
metric appearing in~(\ref{eq0311-1}).
For ease of notation, we denote by $\eps$ and~$h$ two positive
parameters which correspond to $h^{({m})}$ and $\eps^{({m})}$ above.
We denote by $\Sig'$ a square matrix with $\Sig'(\Sig')^*=(\int
_{B(0,h)} x_i x_j \nu(\dd x))_{i,j\in\{1,\ldots,d_X\}}$. Moreover,
we let $L'$ denote the process constituted by the compensated jumps of
$L$ of size larger than $h$, and let $B=(B_t)_{t\ge0}$ be a
$d_X$-dimensional Wiener process that is independent of $W$ and $L'$.
Then we consider the solution $\ups=(\ups_t)_{t\ge0}$ of the
integral equation
\[
\ups_t=y_0+\int_0^t a\bigl(\ups_{\io(s-)}\bigr)\,\dd\cX_s,
\]
where $\cX=(\cX_t)_{t\ge0}$ is given as $\cX_t=\Sig W_t+\Sig'
B_t+L'_t+bt$ and $\io(t)=\max(\II\cap[0,t])$, where $\II$ is, in
analogy to above, the set of random times $(T'_j)_{j\in\IZ_+}$
defined inductively via $T'_0=0$ and
\[
T'_{j+1}=\inf\{ t\in(T'_j,\infty)\dvtx  |\Delta L_t|\ge h\mbox{ or
}t=T'_j+\eps\} \qquad  \mbox{for }j\in\IZ_+.
\]

The process $\ups$ is closely related to $\ups^{({m})}$ from
Section~\ref{sec:alg} and choosing $\eps=\eps_m$ and $h=h_m$,
implies that $(\ups_{\io(t)})_{t\ge0}$ and $\ups^{({m})}$ are
identically distributed.

We need to introduce two further crucial quantities: for $h>0$, let
$F(h)=\int_{B(0,h)} |x|^2 \nu(\dd x)$ and $F_0(h)=\int_{B(0,h)^c}
x \nu(\dd x)$.

\begin{theorem}\label{th0310-1} Suppose that Assumption  \ref{assu1}
is valid. There exists a finite constant $\kappa$ that depends only on
$K$, $d_X$  and $\vth$ such that for $\eps\in(0,\frac12]$, $\eps
'\in[2\eps,1]$, and $h\in(0,\mathfrak h]$ with $\nu(B(0,h)^c) \le
\frac1\eps$, one has
\[
\cW\bigl(Y, \ups_{\io(\cdot)}\bigr)^2\le\kappa \biggl[{F(h) \eps'} +\frac
{h^2}{{\eps'}} \log \biggl(\frac{\eps' F(h)}{h^2}\vee e \biggr)^2+
{\eps\log\frac e\eps} \biggr],
\]
and, if $\Sigma=0$, one has
\[
\cW\bigl(Y, \ups_{\io(\cdot)}\bigr)^2\le\kappa \biggl[F(h)  \biggl(\eps
'+\eps\log\frac e\eps \biggr) +\frac{h^2}{{\eps'}} \log
\biggl(\frac{\eps' F(h)}{h^2}\vee e \biggr)^2+  |b-F_0(h) |^2
\eps^2 \biggr].
\]
\end{theorem}

\begin{cor}\label{cor0311-1} Under Assumption \ref{assu1},
there exists a constant $\kappa=\kappa(K,d_X,\vth)$ such that for
all $\eps\in(0, \frac14]$ and $h\in(0, \mathfrak{h}]$ with $\nu
(B(0,h)^c)\vee\frac{F(h)}{h^2} \le\frac1\eps$, one has
\[
\cW\bigl(Y,\ups_{\io(\cdot)}\bigr)^2\le\kappa \biggl(h^2 \frac1{\sqrt\eps
}+\eps \biggr) \log\frac e\eps,
\]
and, in the case where $\Sig=0$,
\[
\cW\bigl(Y,\ups_{\io(\cdot)}\bigr)^2\le\kappa \biggl(h^2 \frac1{\sqrt\eps
} \log\frac e\eps+  |b-F_0(h) |^2 \eps^2 \biggr).
\]
\end{cor}

\begin{pf}
Choose $\eps'= \sqrt\eps\log1/\eps$ and observe that $\eps'\ge
2\eps$ since $\eps\le\frac14$. Using that $\frac{F(h)}{h^2} \le
g(h)\le\frac1{\eps}$, it is straight forward to verify the estimate
with Theorem~\ref{th0310-1}.
\end{pf}

\subsection[Strategy of the proof of Theorem 3.1 and main
notation]{Strategy of the proof of Theorem~\protect\ref{th0310-1} and main
notation}\label{sec31}

We represent $X$ as
\[
X_t=\Sig W_t+ L'_t+L''_t+bt,
\]
where $L''=(L''_t)_{t\ge0}=L-L'$ is the process which comprises the
compensated jumps of $L$ of size smaller than $h$.
Based on an additional parameter $\eps'\in[2\eps,1]$, we couple
$L''$ with $\Sig B$. The introduction of the explicit coupling is
deferred to Section~\ref{sec:gauss}. Let us roughly explain the idea
behind the parameter $\eps'$.
In classical Euler schemes, the coefficients of the SDE are updated in
either a deterministic or a random number of steps of a given (typical)
length. Our approximation updates the coefficients at steps of order
$\eps$ as the classical Euler method. However, in our case the L\'evy
process that comprises the small jumps is ignored for most of the time
steps. It is only considered on steps of order of size $\eps'$.

On the one hand, a large $\eps'$ reduces the accuracy of the
approximation. On the other hand, the part of the small jumps has to be
approximated by a Wiener process and the error inferred from the
coupling decreases in $\eps'$. This explains the increasing and
decreasing terms in Theorem~\ref{th0310-1}. Balancing $\eps'$ and
$\eps$ then leads to Corollary~\ref{cor0311-1}.

We need some auxiliary processes.
Analogously to $\II$ and $\io$, we let $\IIJ$ denote the set of
random times $(T_j)_{j\in\IZ_+}$ defined inductively by $T_0=0$ and
\[
T_{j+1}=\min \bigl( \II\cap(T_j+\eps'-\eps,\infty) \bigr)
\]
so that the mesh-size of $\IIJ$ is less than or equal to $\eps'$.
Moreover, we set $\eta(t)=\max(\IIJ\cap[0,t])$.

Let us now introduce the first auxiliary processes.
We set $X'=(X_t-L''_t)_{t\ge0}$ and we consider the solution $\x'=(\x
'_t)_{t\ge0}$ to the integral equation
%
\begin{equation}\label{sde2}
\x'_t=y_0+\int_0^t a\bigl(\xism'\bigr)\,\dd X'_s +\int_0^ta\bigl(\xesm'\bigr)\,\dd
L''_{\eta(s)}
\end{equation}
and the process $\x=(\x_t)_{t\ge0}$ given by
\[
\x_t=\x'_t+a\bigl(\xet'\bigr)\bigl(L''_t-L''_{\eta(t)}\bigr).
\]
It coincides with $\x'$ for all times in $\IIJ$ and satisfies
\[
\x_t= y_0+\int_0^t a\bigl(\xism'\bigr)\,\dd X'_s +\int_0^ta\bigl(\xesm\bigr)\,\dd L''_s.
\]

Next, we replace the term $L''$ by the Gaussian term $\Sig' B$ in the
above integral equations and obtain analogs of $\x'$ and $\x$ which
are denoted by $\bar\ups'$ and $\bar\ups$. To be more precise,
$\ups'=(\bar\ups'_t)_{t\ge0}$ is the solution to the stochastic
integral equation
\[
\bar\Upsilon'_t=y_0+\int_0^t a\bigl(\bar\ups'_{\io(s-)}\bigr)\,\dd X'_s
+\int_0^ta\bigl(\bar\ups'_{\eta(s-)}\bigr)\Sigma' \,\dd B_{\eta(s)},
\]
and $\bar\ups=(\bar\ups_t)_{t\ge0}$ is given via
\[
\bar\ups_t= \bar\ups'_t+ a\bigl(\bar\ups'_{\eta(t)}\bigr) \Sigma' \bigl(B_{t}
-B_{\eta(t)}\bigr).
\]
%

We now focus on the discrepancy of $Y$ and $\ups_{\io(\cdot)}$. By
the triangle inequality, one has
%
\begin{equation}\label{eq0218-3}
\bigl\|Y-\ups_{\io(\cdot)}\bigr\|\leq\|Y-\bar Y\| + \|\bar Y- \bar\ups\| +
\|\bar\ups-\ups\|+\bigl\|\ups-\ups_{\io(\cdot)}\bigr\|.
\end{equation}
Moreover, the second term on the right satisfies
%
\begin{equation}\label{eq0218-4}
\|\bar Y- \bar\ups\| \le\|\bar Y'-\bar\ups'\| +\|\bar Y-\bar Y'
-(\bar\ups-\bar\ups')\|.
\end{equation}

In order to prove Theorem~\ref{th0310-1}, we control the error terms
individually. The first term on the right-hand side of (\ref
{eq0218-3}) is considered in Proposition~\ref{propSA1}. The third and
fourth term are treated in Propositions~\ref{prop0306-1} and~\ref
{prop0306-2}, respectively.
The terms on the right-hand side of (\ref{eq0218-4}) are investigated
in Propositions~\ref{prop0218-1} and \ref{prop0218-2}, respectively.
Note that only the latter two expressions depend on the particular
choice of the coupling of $L''$ and $\Sig' B$. Once the
above-mentioned propositions are proved, the statement of Theorem~\ref
{th0310-1} follows immediately by combining these estimates and
identifying the dominant terms.

\section{Approximation of $Y$ by $\bar Y$}\label{sec:4}

\begin{prop} \label{propSA1} There exists a constant $\kappa>0$
depending on $K$ only such that, for $\eps\in(0,\frac12]$, $\eps
'\in[2\eps,1]$  and $h>0$ with $\nu(B(0,h)^c)\le\frac1\eps$, one has
\[
\IE \Bigl[\sup_{t\in[0,1]}|Y_t-\x_t|^2 \Bigr] \le\kappa
[ F(h)  \eps' + |b-F_0(h) |^2\eps^2  ],
\]
if $\Sigma=0$, and
%
\begin{equation}\label{sde3}
\IE \Bigl[\sup_{t\in[0,1]}|Y_t-\x_t|^2 \Bigr] \leq\kappa \bigl(
\eps+ F(h) \eps' \bigr)
\end{equation}
for general $\Sigma$.
\end{prop}

\begin{pf}
For $t\ge0$, we consider $Z_t=Y_t-\x_t$, $Z'_t=Y_t-\xit'$,
$Z''_t=Y_t-\xet$  and $z(t)=\IE[\sup_{s\in[0,t]} |Z_s|^2]$. The
main task of the proof is to establish an estimate of the form
\[
z(t)\leq\alpha_1 \int_0^t z(s)\,\dd s + \alpha_2
\]
for appropriate values $\alpha_1,\alpha_2>0$. Since $z$ is finite
(see, for instance, \cite{DerHei10}), then Gronwall's inequality
implies as upper bound:
\[
\IE\Bigl[\sup_{s\in[0,1]} |Y_s-\x_s|^2\Bigr] \leq\alpha_2  \exp({\alpha_1}).
\]
We proceed in two steps.

\textit{1st step.}
Note that
\begin{eqnarray*}
Z_t&=& \underbrace{ \int_0^t \bigl(a(Y_{s-})-a\bigl(\xism'\bigr)\bigr)\, \dd(\Sigma
W_s+L'_s) +\int_0^t\bigl(a(Y_{s-})-a\bigl(\xesm\bigr)\bigr)\, \dd L''_s}_{=:M_t} \\
&&{}    +\int_0^t \bigl(a(Y_{s-})-a\bigl(\xism'\bigr)\bigr)b\,\dd s,
\end{eqnarray*}
so that
%
\begin{equation}\label{eq0304-1}
|Z_t|^2 \le2   |M_t|^2 +2    \biggl| \int_0^t \bigl(a(Y_{s-})-a\bigl(\bar
Y_{\io(s-)} \bigr)\bigr)b\,\dd s \biggr|^2.
\end{equation}
For $t\in[0,1]$, we conclude with the Cauchy--Schwarz inequality that
the second term on the right-hand side is bounded by
$2 K^4 \int_0^t |Z'_{s-}|^2\,\dd s$.

Certainly, $(M_t)$ is a (local) martingale with respect to the
canonical filtration, and we apply the Doob inequality together with
Lemma~\ref{le0713-1} to deduce that
\begin{eqnarray*}
\IE \Bigl[\sup_{s\in[0,t]} |M_s|^2 \Bigr]&\le&4  \IE\biggl [ \int
_0^t \bigl|a(Y_{s-})-a\bigl(\xism'\bigr)\bigr|^2\,\dd
\langle\Sigma W+L'\rangle_s\\
&&\hspace*{32pt}{}+\int
_0^t \bigl|a(Y_{s-})-a\bigl(\xesm\bigr)\bigr|^2\,\dd\langle L''\rangle_s \biggr].
\end{eqnarray*}
Here and elsewhere, for a multivariate local $L^2$-martingale
$S=(S_t)_{t\ge0}$, we denote $\langle S\rangle=\sum_j \langle
S^{({j})}\rangle$ and $\langle S^{({j})}\rangle$ denotes the
predictable compensator of the classical bracket process of the $j$th
coordinate $S^{({j})}$ of $S$.
Note that $\dd\langle\Sigma W+L'\rangle_t=(|\Sig|^2+\int
_{B(0,h)^c}|x|^2 \nu(\dd x))\, \dd t\le2 K^2 \,\dd t$ and $\dd
\langle L''\rangle_t=F(h) \, \dd t$. Consequently,
\[
\IE \Bigl[\sup_{s\in[0,t]} |M_s|^2 \Bigr]\le4   \IE \biggl[2K^4
\int_0^t |Z'_s|^2\,\dd s+K^2  F(h) \int_0^t |Z''_s|^2\,\dd s \biggr].
\]
Hence, by (\ref{eq0304-1}) and Fubini's theorem, one has
\[
\IE \Bigl[\sup_{s\in[0,t]} |Z_s|^2  \Bigr] \leq\kappa_1 \int_0^t
 \bigl[z(s) + \IE[|Z'_{s}|^2 ]+F(h)  \IE[|Z''_s|^2 ]\bigr]\,\dd s
\]
for a constant $\kappa_1$ that depends only on $K$.
Since $Z'_t=Z_t+\bar Y_t-\bar Y'_{\io(t)}$ and $Z''_t=Z_t+\bar
Y_t-\bar Y_{\eta(t)}$, we get
%
\begin{equation}\label{eq0612-1}
z(t)\leq\kappa_2 \int_0^t  \bigl[z(s) + \IE\bigl[\bigl|\x_s-\xis'\bigr|^2\bigr] +
F(h)  \IE\bigl[\bigl|\x_s-\xes\bigr|^2\bigr] \bigr]\,\dd s
\end{equation}
for an appropriate constant $\kappa_2=\kappa_2(K)$.

\textit{2nd step.} In the second step we provide appropriate estimates
for $\IE[|\x_t-\xit'|^2]$ and $\IE[|\x_t-\xet|^2]$. The processes
$W$ and $L''$ are independent of the random time $\io(t)$. Moreover,
$L'$ has no jumps in $(\io(t),t)$, and we obtain
\begin{eqnarray*}
\x_t- \xit' &=&\x'_t-\xit'+a\bigl(\xet\bigr)\bigl(L''_{t}-L''_{\eta(t)}\bigr) \\
&=& a\bigl(\xit'\bigr)\bigl ( \Sigma\bigl(W_t-W_{\io(t)}\bigr) +  \bigl(b-F_0(h) \bigr)
\bigl(t-\io(t)\bigr) \bigr)\\
&&{}    +a\bigl(\xet\bigr) \bigl(L''_t-L''_{\eta(t)}\bigr)
\end{eqnarray*}
so that
\begin{eqnarray*}
\IE\bigl[\bigl|\x_t- \xit'\bigr|^2\bigr] &\le&3 K^2  \bigl[  \IE\bigl[\bigl(\bigl|\bar Y'_{\io
(t)}-y_0\bigr|+1\bigr)^2\bigr]  \bigl( |\Sig|^2\eps+ |b-F_0(h) |^2 \eps
^2 \bigr)\\
&&\hspace*{89pt}{}   +\IE\bigl[\bigl(\bigl|\bar Y_{\eta(t)}-y_0\bigr|+1\bigr)^2\bigr] F(h)  \eps'   \bigr].
\end{eqnarray*}
By Lemma~\ref{le0305-1}, there exists a constant $\kappa_3=\kappa
_3(K)$ such that
%
\begin{equation}\label{eq0305-1}
\IE\bigl[\bigl|\x_t- \xit'\bigr|^2\bigr] \le\kappa_3 [ |\Sig|^2\eps+
|b-F_0(h) |^2 \eps^2+ F(h)  \eps' ].
\end{equation}

Similarly, we estimate $\IE[|\x_t-\xet|^2]$.
Given $\eta(t)$, $(L'_{\eta(t)+u} -\break L'_{\eta(t)})_{u\in[0,(\eps
'-\eps) \wedge  (t-\eta(t))]}$ is distributed as the unconditioned L\'
evy process $L'$ on the time interval $[0, (\eps'-\eps) \wedge
(t-\eta(t))]$. Moreover, we have $\dd L'_{u} = -F_0(h) \,\dd u$ on
$(\eta(t)+ \eps'-\eps,t]$. Consequently,
\begin{eqnarray*}
\x_t-\xet&=& \int_{\eta(t)}^t \ind_{\{s-\eta(t)\le\eps'-\eps\}}
a\bigl(\xism\bigr)\,\dd (\Sigma W_s+L'_s+b s )\\
&&{}    + \int_{\eta(t)}^t \ind_{\{s-\eta(t)> \eps'-\eps\}}
a\bigl(\xism\bigr)\,\dd \bigl(\Sigma W_s+\bigl(b-F_0(h)\bigr)s \bigr)\\
&&{}    + a\bigl(\xet\bigr)\bigl (L''_t-L''_{\eta(t)}\bigr),
\end{eqnarray*}
and analogously as we obtained (\ref{eq0305-1}) we get now that
\[
\IE\bigl[\bigl|\x_t-\xet\bigr|^2\bigr]\le\kappa_4 [ \eps' +  |b-F_0(h)
|^2 \eps^2 ]
\]
for a constant $\kappa_4=\kappa_4(K)$. Next, note that, by the
Cauchy--Schwarz inequality, $|F_0(h)|^2\leq\int_{B(0,h)^c} |x|^2 \nu
(\dd x)\cdot\nu(B(0,h)^c) \leq\frac{K^2}{\eps}$ so that we arrive at
\[
\IE\bigl[\bigl|\x_t-\xet\bigr|^2\bigr]\le\kappa_5   \eps'.
\]

Combining this estimate with (\ref{eq0612-1}) and (\ref{eq0305-1}),
we obtain
\[
z(t)\le\kappa_2 \int_0^t z(s)\,\dd s +\kappa_6  [
|\Sigma|^2 \eps+ F(h) \eps' +|b-F_0(h) |^2\eps^2  ].
\]
In the case where $\Sigma=0$, the statement of the proposition follows
immediately via Gronwall's inequality.
For general $\Sigma$, we obtain the result by recalling that
$|F_0(h)|^2\leq\frac{K^2}{\eps}$.
\end{pf}

\section{Approximation of $\bar\ups$ by $\ups_{\io(\cdot)}$}
\label{sec:5}

\begin{prop}\label{prop0306-1} Under the assumptions of
Proposition~\ref{propSA1}, one has
\[
\IE[\|\bar\ups-\ups\|^2]\le\kappa\eps' F(h)
\]
for a constant $\kappa$ depending only on $K$.
\end{prop}

\begin{pf}
The proposition can be proved as Proposition~\ref{propSA1}. Therefore,
we only provide a sketch of the proof.
The arguments from the first step give, for $t\in[0,1]$,
\[
z(t)\le\kappa_1 \int_0^t \bigl [ z(s) + \IE\bigl[\bigl|\bar\ups_{\io
(s)}-\bar\ups'_{\io(s)}\bigr|^2\bigr] + F(h)   \IE\bigl[\bigl|\bar\ups_{\io(s)}
-\bar\ups_{\eta(s)}\bigr|^2\bigr] \bigr]\,\dd s,
\]
where $z(t)=\IE[\sup_{s\in[0,t]}|\ups_s-\bar\ups_{s}|^2]$ and
$\kappa_1=\kappa_1(K)$ is an appropriate constant.

Moreover, based on Lemma~\ref{le0305-1} the second step leads to
\[
\IE\bigl[\bigl|\bar\ups_{\io(t)}-\bar\ups'_{\io(t)}\bigr|^2\bigr]\le\kappa_2 \eps
' F(h)  \quad \mbox{and}\quad   \IE\bigl[\bigl|\bar\ups_{\io(t)}-\bar\ups_{\eta
(t)}\bigr|^2\bigr]\le\kappa_3 \eps'
\]
for appropriate constants $\kappa_2=\kappa_2(K)$ and $\kappa
_3=\kappa_3(K)$.
Then Gronwall's lemma implies again the statement of the proposition.
\end{pf}

\begin{prop}\label{prop0306-2} Under the assumptions of
Proposition~\ref{propSA1}, there
exists a constant $\kappa$ depending only on $K$ and $d_X$ such that,
if $\Sig=0$,
\[
\IE\Bigl[ \sup_{t\in[0,1]} \bigl|\ups_t-\ups_{\io(t)}\bigr|^2\Bigr]\le\kappa \biggl[
F(h)  \eps\log\frac e\eps+  |b-F_0(h) |^2 \eps^2 \biggr]
\]
and, in the general case,
\[
\IE\Bigl[ \sup_{t\in[0,1]} \bigl|\ups_t-\ups_{\io(t)}\bigr|^2\Bigr]\le\kappa
\eps\log\frac e\eps.
\]
\end{prop}

\begin{pf} Recall that by definition
\[
\ups_t-\ups_{\io(t)}=\int_{\io(t)}^t a\bigl(\ups_{\io(s-)}\bigr) \,\dd\cX_s
\]
so that
\[
\bigl|\ups_t-\ups_{\io(t)}\bigr|^2\le K^2 \bigl(\bigl|\ups_{\io(t)}-y_0\bigr|+1\bigr)^2   \bigl|\cX
_t-\cX_{\io(t)}\bigr|^2.
\]
Next, we apply Lemma~\ref{le1121-1}. For $j\in\IZ_+$, we choose
\[
U_j= |\ups_{T'_j\wedge1}-y_0|^2 \quad  \mbox{and}\quad   V_j= \sup_{s\in
[T'_{j} ,T'_{j+1}\wedge1)} \bigl|\cX_{t}-\cX_{\io(t)}\bigr|^2
\]
with the convention that the supremum of the empty set is zero. Then
\begin{eqnarray*}
\IE\Bigl[\sup_{t\in[0,1]} \bigl|\ups_t-\ups_{\io(t)}\bigr|^2\Bigr]
&\le&\IE
\Bigl[\sup_{j\in\IZ_+} U_j \Bigr] \cdot\IE \Bigl[\sup_{j\in\IZ_+}
V_j \Bigr] \\
&\le&\IE \Bigl[\sup_{t\in[0,1]}(|\ups_{t}-y_0|+1)^2 \Bigr] \cdot
\IE \Bigl[\suptwo{0\le s<t\le1}{t-s\le\eps} |\cX_{t}-\cX_s|^2
 \Bigr].
\end{eqnarray*}
By Proposition~\ref{prop0306-1} and Lemma~\ref{le0305-1}, $\IE[\sup
_{t\in[0,1]} (|\ups_{t}-y_0|+1)^2 ]$ is bounded by a constant that
depends only on $K$.

Consider $\vphi\dvtx [0,1]\to[0,\infty), \delta\mapsto\sqrt{\delta
\log(e/\delta)}$. By L\'evy's modulus of continuity,
\[
\|W\|_{\vphi}:=\sup_{0\le s<t\le1} \frac{|W_t-W_s|}{\vphi(t-s)}
\]
is finite almost surely, so that Fernique's theorem implies that $\IE
[\|W\|_\vphi^2]$ is finite too.
Consequently,
%
\begin{eqnarray}\label{eq0512-1}
&&\IE\Bigl[ \sup_{s\in[0,t]} \bigl|\cX_s-\cX_{\io(s)}\bigr|^2\Bigr]\nonumber
\\[-8pt]
\\[-8pt]
&&\qquad\le3  \biggl[ \bigl(|\Sig
|^2+F(h)\bigr) \IE[\|W\|_\vphi^2]  \eps\log\frac e\eps+
|b-F_0(h) | ^2 \eps^2 \biggr].\nonumber
\end{eqnarray}
The result follows immediately by using that $|F_0(h)|^2\le\frac
{K^2}{\eps}$ and ruling out the asymptotically negligible terms.
\end{pf}

\section{Gaussian approximation via Koml\'{o}s, Major  and Tusn\'{a}dy}
\label{sec:KMT}

In this section, we prove
the following theorem.
\begin{theorem} \label{th0217-1} Let $h>0$ and $L=(L_t)_{t\ge0}$ be a
$d$-dimensional $(\nu,0)$-L\'evy martingale whose L\'evy measure $\nu
$ is supported on $B(0,h)$. Moreover, we suppose that
for $\vth\ge1$, one has
\[
\int\langle y',x\rangle^2  \nu(\dd x) \le\vth\int\langle
y,x\rangle^2  \nu(\dd x)
\]
for any $y,y'\in\IR^d$ with $|y|=|y'|$, and set $\sigma^2=\int
|x|^2 \nu(\dd x)$.

There exist constants $c_1,c_2>0$ depending only on $d$ such that the
following statement is true.
For every $T\ge0$, one can couple the process $(L_t)_{t\in[0,T]}$
with a Wiener process $(B_t)_{t\in[0,T]}$ such that
\[
\IE\exp \biggl\{\frac{c_1}{\sqrt{\vth}h } \sup_{t\in[0,T]}
|L_t-\Sigma B_t|  \biggr\} \le\exp \biggl\{c_2 \log \biggl( \frac{
\sigma^2 T}{h^2}\vee e \biggr) \biggr\},
\]
where $\Sigma$ is a square matrix with $\Sigma\Sigma^*=\cov_{L_1}$
and $\sigma^2=\int|x|^2 \nu(\dd x)$.
\end{theorem}

The proof of the theorem is based on Zaitsev's generalization \cite
{Zai98} of the Koml\'os--Major--Tusn\'ady coupling. In this context, a
key quantity is the \textit{Zaitsev parameter}: Let $Z$ be a
$d$-dimensional random variable with finite exponential moments in a
neighborhood of zero and set
\[
\Lam(\th)=\log\IE\exp\{\langle\th,Z\rangle\}
\]
for all $\th\in\IC$ with integrable expectation.
Then the parameter is defined as
\begin{eqnarray*}
\tau(Z)= \inf \{ \tau>0\dvtx  |\partial_w\, \partial_v^2 \Lam(\th
)|&\le&\tau \langle\cov_Z  v,v\rangle\mbox{ for all }\th\in\IC
^d, v,w\in\IR^d \\
&&\hspace*{21.5pt}\mbox{with }|\th|\le\tau^{-1}\mbox{ and }
|w|=|v|=1 \}.
\end{eqnarray*}
In the latter set, we implicitly only consider $\tau$'s for which
$\Lam$ is finite on a neighborhood of $\{x\in\IC^d\dvtx  |x|\le1/\tau\}
$. Moreover, $\cov_Z$ denotes the covariance matrix of~$Z$.

\begin{pf*}{Proof of Theorem~\ref{th0217-1}}
\textit{1st step:} First, consider a $d$-dimensional infinitely
divisible random variable $Z$ with
\[
\Lam(\th):=\log\IE e^{\langle\th,Z\rangle}= \int\bigl(e^{\langle\th
,x\rangle}-\langle\th,x\rangle-1\bigr) \nu'(\dd x),
\]
where the L\'evy measure $\nu'$ is supported on the ball $B(0,h')$ for
a fixed $h'>0$.
Then
\[
\partial_w \,\partial_v^2\Lam(\th)= \int_{B(0,h')} \langle
w,x\rangle\langle v,x\rangle^2 e^{\langle\th, x\rangle} \nu(\dd x)
\]
and
\[
\langle\cov_{Z}  v,v\rangle=\var\langle v, Z\rangle= \partial
_v^2 \Lam_{Z}(0)=\int_{B(0,h')} \langle v,x\rangle^2  \nu(\dd x).
\]
We choose $\zeta>0$ with $e^{\zeta}=1/\zeta$, and observe that for
any $\th\in\IC^d, v,w\in\IR^d$ with $|\th|\le\zeta/ h'$ and $|w|=|v|=1$,
\[
|\partial_w \,\partial_v^2\Lam(\th)|\le h' e^{|\th| h'} \langle\cov
_{Z}  v,v\rangle\le\frac{h'}\zeta  \langle\cov_{Z}  v,v\rangle.
\]
Hence,
\[
\tau(Z) \le\frac{h'}\zeta.
\]

\textit{2nd step:} In the next step, we apply Zaitsev's coupling to
piecewise constant interpolations of $(L_t)$. Fix $m\in\IN$ and
consider $L^{({m})}=(L^{({m})}_t)_{t\in[0,T]}$ given via
\[
L^{({m})}_t= L_{\lfl2^m t/T\rfl2^{-m}T}.
\]
Moreover, we consider a $d$-dimensional Wiener process $B=(B_t)_{t\ge
0}$ and its piecewise constant interpolation $\Sigma B^{({m})}$ given
by $B^{({m})}=(B_{\lfl2^m t/T\rfl2^{-m}T})_{t\in[0,T]}$.

Since $\cov_{L_1}$ is self-adjoint, we find a representation $\cov
_{L_t}=t U D U^*$ with $D$ diagonal and $U$ orthogonal. Hence, for
$A_t:=(tD)^{-1/2} U^*$ we get
$\cov_{A_t L_t}=I_d$. We denote by $\lambda_1$ the leading and by
$\lambda_2$ the minimal eigenvalue of $D$ (or $\cov_{L_1}$). Then
$A_t L_t$ is again infinitely divisible and the corresponding L\'evy
measure is supported on $B(0,h/\sqrt{\lambda_2 t})$. By part one, we
conclude that
\[
\tau(A_t L_t)\le\frac{h}{\zeta\sqrt{\lambda_2 t}}.
\]

Now the discontinuities of $A_{2^{-m}} L^{({m})}$ are i.i.d. with
unit covariance and Zaitsev parameter less than or equal to $\frac
{h2^{m/2}}{\zeta\sqrt{T\lambda_2}}$. By \cite{Zai98}, Theorem~1.3,
one can couple the processes $L$ and $\Sigma B$ on an appropriate
probability space such that
\[
\IE\exp \biggl\{ \kappa_1 \frac{\sqrt{T\lambda_2}}{2^{m/2} h} \sup
_{t\in[0,T]} \bigl|A_{2^{-m}}L^{({m})}_t-A_{2^{-m}}\Sigma B^{({m})}_t\bigr|
\biggr\}\le\exp\biggl \{ \kappa_2 \log \biggl(\frac{\zeta^2
T\lambda_2}{h^2}\vee e \biggr) \biggr\},
\]
where $\kappa_1,\kappa_2>0$ are constants only depending on the
dimension $d$.
The smallest eigenvalue of $A_{2^{-m}}$ is $2^{m/2} (T\lambda
_1)^{-1/2}$ and, by assumption, $\lambda_1\le\vth\lambda_2$. Since
$\lambda_2\le\sigma^2$, we get
\[
\IE\exp \biggl\{ \kappa_1 \frac1{\sqrt{\vth}h} \sup_{t\in[0,T]}
\bigl|L^{({m})}_t-\Sigma B^{({m})}_t\bigr| \biggr\}\le\exp \biggl\{ \kappa_2
\log \biggl(\frac{\zeta^2 T\sigma^2}{h^2}\vee e \biggr) \biggr\}.
\]

\textit{3rd step:}
The general result follows by approximation. First, note that\break  $\sup
_{t\in[0,T]} |L_t-L^{({m})}_t|$ converges as $m\to\infty$ to $\sup
_{t\in[0,T]} |L_t-L_{t-}|$ so that by dominated convergence
\begin{eqnarray*}
&&\lim_{m\to\infty} \IE\exp \biggl\{ \kappa_1 \frac1{\sqrt{\vth
}h} \sup_{t\in[0,T]} \bigl|L_t-L^{({m})}_t\bigr| \biggr\}\\
&&\qquad= \IE\exp\biggl \{
\kappa_1 \frac1{\sqrt{\vth}h} \sup_{t\in[0,T]} |L_t-L_{t-}|
\biggr\}\le e^{\kappa_1}.
\end{eqnarray*}
Analogously, $\lim_{m\to\infty} \IE\exp \{ \kappa_1 \frac
1{\sqrt{\vth}h} \sup_{t\in[0,T]} |\Sigma B_t-\Sigma
B^{({m})}_t| \}=1$.
Next, we choose $\kappa_3\ge1$ with $e^{\kappa_1}+1\le e^{\kappa
_2+\kappa_3}$ and we fix $m\in\IN$ such that
\begin{eqnarray*}
&&\IE\exp \biggl\{ \sfrac{\kappa_1}3 \frac1{\sqrt{\vth}h} \sup
_{t\in[0,T]} |L_t-\Sigma B_t| \biggr\} + \IE\exp \biggl\{ \kappa_1
\frac1{\sqrt{\vth}h} \sup_{t\in[0,T]} \bigl|\Sigma B_t-\Sigma
B^{({m})}_t\bigr| \biggr\}\\
&&\qquad \le e^{\kappa_2+\kappa_3}.
\end{eqnarray*}
We apply the coupling introduced in step~2 and estimate
\begin{eqnarray*}
\IE\exp \biggl\{ \sfrac{\kappa_1}3 \frac1{\sqrt{\vth}h} \sup
_{t\in[0,T]} |L_t-\Sigma B_t| \biggr\}
&\le&\IE\exp \biggl\{ \kappa_1
\frac1{\sqrt{\vth}h} \sup_{t\in[0,T]} \bigl|L_t-L^{({m})}_t\bigr| \biggr\}\\
&&{}    +\IE\exp \biggl\{ \kappa_1 \frac1{\sqrt{\vth}h} \sup_{t\in
[0,T]} \bigl|L_t^{({m})}-\Sigma B^{({m})}_t\bigr| \biggr\}\\
&&{}   +\IE\exp \biggl\{ \kappa_1 \frac1{\sqrt{\vth}h} \sup_{t\in
[0,T]} \bigl|\Sigma B^{({m})}_t-\Sigma B_t\bigr| \biggr\}\\
&\le&\exp\biggl \{ \kappa_2 \log \biggl(\frac{T\sigma^2}{h^2}\vee
e \biggr) \biggr\} + e^{\kappa_2+\kappa_3}.
\end{eqnarray*}
Straightforwardly, one obtains the assertion of the theorem for
$c_1=\kappa_1/3$ and $c_2=\kappa_2+2\kappa_3$.
\end{pf*}

\begin{cor}\label{cor0306-1} The coupling introduced in Theorem~\ref
{th0217-1} satisfies
\[
\IE\Bigl[\sup_{t\in[0,T]} |L_t -\Sig B_t |^2\Bigr]^{1/2}\le\frac{\sqrt\vth
h}{c_1} \biggl(c_2 \log\biggl (\frac{\sigma^2 T}{h^2} \vee e
\biggr)+2 \biggr),
\]
where $c_1$ and $c_2$ are as in the theorem.
\end{cor}

\begin{pf}
We set $Z=\sup_{t\in[0,T]} |L_t -\Sig B_t |$ and $t_0=\frac{\sqrt
\vth h}{c_1}c_2 \log (\frac{\sigma^2T}{h^2}\vee e )$, and
use that
%
\begin{equation}\label{eq0212-1}
\IE[ Z^2]= 2 \int_0^\infty t \IP(Z\ge t)\,\dd t \le t_0^2 + 2 \int
_{t_0}^\infty t \IP(Z\ge t)\,\dd t.
\end{equation}
By the Markov inequality and Theorem~\ref{th0217-1}, one has for $s\ge0$
\[
\IP(Z\ge s+t_0)\le\frac{\IE [\exp \{ {c_1}/({\sqrt
{\vth}h }) Z  \} ]}{\exp \{ {c_1}/({\sqrt{\vth}h })
(s+t_0)  \}}\le\exp \biggl\{-\frac{c_1}{\sqrt\vth h} s \biggr\}.
\]
We set $\alpha=\sqrt\vth h/c_1$, and deduce together with (\ref
{eq0212-1}) that
\[
\IE[ Z^2]\le t_0^2+ 2 \int_{0}^\infty(s+t_0) \exp\biggl\{-\sfrac1\alpha
s\biggr\}\,\dd s = t_0^2 +2 t_0 \alpha+2 {\alpha^2}\le(t_0+2 \alpha)^2.\quad
\]
\upqed
\end{pf}

\section{Coupling the Gaussian approximation}\label{sec:gauss}

We are now in the position to couple the processes $L''$ and $\Sig' B$
introduced in Section~\ref{sec31}. We adopt again the notation of
Section~\ref{sec31}.

To introduce the coupling, we need to assume that Assumption \ref{assu1}
is valid, and that $\eps\in(0,\frac12]$, $\eps'\in[2\eps,1]$  and
$h\in(0,\mathfrak h]$ are such that $\nu(B(0,h)^c)\le\frac1\eps$.
Recall that $L''$ is independent of $W$ and $L'$.
In particular, it is independent of the times in $\IIJ$, and given $W$
and $L'$ we couple the Wiener process $B$ with $L''$ on each interval
$[T_i,T_{i+1}]$ according to the coupling provided by Theorem~\ref{th0217-1}.

More explicitly, the coupling is established in such a way that, given
$\IIJ$, each pair of processes $(B_{t+T_j}-B_{T_j})_{t\in
[0,T_{j+1}-T_j]}$ and $(L''_{t+T_j}-L''_{T_j})_{t\in[0,T_{j+1}-T_j]}$
is independent of $W$, $L'$  and the other pairings, and satisfies
%
\begin{eqnarray}\label{eq0325-1}
&&\IE \biggl[\exp \biggl\{\frac{c_1}{\sqrt{\vth}h } \sup_{t\in
[T_j,T_{j+1}]} |L''_t-L''_{T_j}-(\Sigma' B_t- \Sigma' B_{T_j}) |
 \biggr\}  \Big| \IIJ \biggr]\nonumber
 \\[-8pt]
 \\[-8pt]
&& \qquad\le\exp \biggl\{c_2 \log \biggl( \frac{
F(h)(T_{j+1}-T_j) }{h^2}\vee e \biggr) \biggr\} \nonumber
\end{eqnarray}
for positive constants $c_1$ and $c_2$ depending only on $d_X$, see
Theorem~\ref{th0217-1}. In particular, by Corollary~\ref{cor0306-1},
one has
%
\begin{eqnarray}\label{eq0218-2}
&&\IE \Bigl[ \sup_{t\in[T_j,T_{j+1}]} |L''_t-L''_{T_j}-(\Sigma' B_t-
\Sigma' B_{T_j}) |^2  | \IIJ \Bigr]^{1/2}\nonumber
\\[-8pt]
\\[-8pt]
&&\qquad\le c_3  h \log \biggl(
\frac{ F(h)(T_{j+1}-T_j) }{h^2}\vee e \biggr)\nonumber
\end{eqnarray}
for a constant $c_3=c_3(d_X,\vth)$.

\begin{prop}\label{prop0218-1}
Under Assumption~\ref{assu1}, there exists a constant $\kappa$ depending only
on $K$, $\vth$  and $d_X$ such that for any $\eps\in(0,\frac12]$,
$\eps'\in[2\eps,1]$  and $h\in(0,\mathfrak{h}]$ with $\nu
(B(0,h)^c) \leq\frac1\eps$, one has
\[
\IE \Bigl[\sup_{[0,1]} |\x'_t-\z'_t|^2 \Bigr]\leq\kappa\frac
1{\eps'}  h^2 \log\biggl (\frac{ \eps' F(h)}{h^2}\vee e \biggr)^2.
\]
\end{prop}


\begin{pf}
For ease of notation, we write
\[
A_t=L''_{\eta(t)}\quad \mbox{and}\quad  A'_t=\Sig' B _{\eta(t)}.
\]
By construction, $(A_t)$ and $(A'_t)$ are martingales with respect to
the filtration $(\cF_t)$ induced by the processes $(W_t)$, $(L'_t)$,
$(A_t)$  and $(A'_t)$.
Let $Z_t=\x'_t-\z'_t$, $Z'_t=\xit'- \zit' $, $Z''_t=\xet'-\zet'$
and $z(t)=\IE[\sup_{s\in[0,t]} |Z_s|^2]$. The proof is similar to
the proof of Proposition~\ref{propSA1}.

Again, we write
%
{\fontsize{10.3}{11.3}\selectfont{\begin{eqnarray}\label{eq0305-5}
Z_t&=&
 \underbrace{ \int_0^t \bigl(a\bigl(\xism'\bigr)-a\bigl(\zism'\bigr)\bigr)\,
 \dd(\Sigma W_s
+L'_s) + \int_0^t a\bigl(\xesm'\bigr)\, \dd A_s- \int_0^t a\bigl(\zesm'\bigr)\,\dd
A'_s}_{=:M_t  \ (\mathrm{local martingale})} \nonumber
\\[-8pt]
\\[-8pt]
&&{}    +\int_0^t \bigl(a\bigl(\xis'\bigr)-a\bigl(\zis'\bigr)\bigr)b \,\dd s.\nonumber
\end{eqnarray}}}
$\!$Denoting $M'= \Sigma W +L'$, we get
\begin{eqnarray*}
 \dd M_t&=& \bigl(a\bigl(\xitm'\bigr)-a\bigl(\zitm'\bigr)\bigr)\, \dd M'_t +a\bigl(\xetm'\bigr)\, \dd
(A_t-A'_t)\\
&& {}+ \bigl(a\bigl(\xetm'\bigr)-a\bigl(\zetm' \bigr)\bigr) \,\dd A'_t
\end{eqnarray*}
and, by Doob's inequality and Lemma~\ref{le0713-1}, we have
%
\begin{eqnarray}\label{eq0305-3}
\IE\Bigl[\sup_{s\in[0,t]} |M_s|^2 \Bigr]
&\le&\kappa_1  \biggl[ \IE \biggl[\int
_0^t |Z'_{s-}|^2 \,\dd\langle M'\rangle_s  \biggr]+ \IE \biggl[ \int
_0^t |Z''_{s-}|^2\, \dd\langle A'\rangle_s \biggr]\nonumber
\\[-8pt]
\\[-8pt]
&&\hspace*{54pt}{} + \IE \biggl[\int
_0^t \bigl(\bigl|\xesm'\bigr|+1\bigr)^2\, \dd\langle A-A'\rangle_s \biggr] \biggr].\nonumber
\end{eqnarray}
Each bracket $\langle\cdot\rangle$ in the latter formula can be
chosen with respect to a (possibly different) filtration such that the
integrand is predictable and the integrator is a local $L^2$-martingale.
As noticed before, with respect to the canonical filtration $(\cF_t)$
one has $\dd\langle M'\rangle_t = (|\Sigma|^2 + \int_{B(0,h)^c}
|x|^2  \nu(\dd x))\,\dd t\le2 K^2 \,\dd t$. Moreover, we have with
respect to the enlarged filtration $(\cF_t \vee\sigma(\IIJ))_{t\ge0}$,
\[
\langle A'\rangle_t=\sum_{\{j\in\IN: T_j\le t\}} (T_j-T_{j-1})
F(h)= \max(\IIJ\cap[0,t]) \cdot F(h),
\]
and, by (\ref{eq0218-2}), for $j\in\IN$,
\[
\Delta\langle A-A'\rangle_{T_j}= \IE [
|L''_{T_j}-L''_{T_{j-1}}-(\Sig' B_{T_j}-\Sig' B_{T_{j-1}})|^2
|   \IIJ ]\le c_3^2 \xi^2,
\]
where $\xi:=h \log(\frac{ \eps' F(h)}{h^2}\vee e)$.
Note that two discontinuities of $\langle A-A'\rangle$ are at least
$\eps'/2$ units apart and the integrands of the last two integrals in
(\ref{eq0305-3}) are constant on $(T_{j-1},T_j]$ so that altogether
\begin{eqnarray*}
\IE\Bigl[\sup_{s\in[0,t]} |M_s|^2 \Bigr] &\le&\kappa_1  \biggl[ 2K^2 \IE
\biggl[ \int_0^t |Z'_s|^2 \,\dd s \biggr] + F(h)  \IE \biggl[\int_0^t
|Z''_s|^2\, \dd s \biggr] \\
&&\hspace*{54pt}{}+ c_3^2 \xi^2 \frac2{\eps'}  \IE \biggl[
\int_0^t \bigl(\bigl|\xesm'\bigr|+1\bigr)^2\, \dd s \biggr] \biggr].
\end{eqnarray*}
With Lemma~\ref{le0305-1} and Fubini's theorem, we arrive at
\[
\IE \Bigl[\sup_{s\in[0,t]} |M_s|^2  \Bigr] \le\kappa_2  \biggl[
\int_0^t z(s) \,\dd s + \xi^2 \frac1{\eps'}  \biggr].
\]
Moreover, by Jensen's inequality, one has
\[
\IE \biggl[\sup_{s\in[0,t]}  \biggl|\int_0^s \bigl(a\bigl(\bar Y'_{\io
(u-)}\bigr)-a\bigl(\bar\ups'_{\io(u-)}\bigr)\bigr)b \,\dd u \biggr|^2 \biggr]\le K^4 \int
_0^t \IE[ |Z'_{s-}|^2]\,\dd s.
\]
Combining the latter two estimates with (\ref{eq0305-5}) and applying
Gronwall's inequality yields the statement of the proposition.
\end{pf}

\begin{prop}\label{prop0218-2} There exists a constant $\kappa$
depending only on $K$ and $d_X$ such that
\begin{eqnarray*}
\IE[\|\bar Y-\bar Y'-(\bar\ups-\bar\ups')\|^2]^{1/2}&\le&\kappa
 \biggl[h\biggl [ \log \biggl(1+\frac2{\eps'} \biggr) + \log
\biggl(\frac{F(h)\eps'}{h^2}\vee e \biggr) \biggr]\nonumber
\\[-8pt]
\\[-8pt]
&&\hspace*{16pt}{}+\sqrt{F(h) \eps'\log
\frac e{\eps'}}  \IE[\|\bar Y'-\bar\ups'\|^2]^{1/2} \biggr].
\end{eqnarray*}
\end{prop}

\begin{pf}
Note that
\begin{eqnarray*}
\bar Y_t-\bar Y'_t-(\bar\ups_t-\bar\ups_t')
&=& a\bigl(\bar Y'_{\eta(t)}\bigr)
\bigl(L_t''-L''_{\eta(t)}\bigr) - a\bigl(\bar\ups'_{\eta(t)}\bigr)\bigl (\Sig' B _t-\Sig'
B _{\eta(t)}\bigr)\\
&=& a\bigl(\bar Y'_{\eta(t)}\bigr)\bigl (L_t''-L''_{\eta(t)}-\bigl(\Sig' B _t-\Sig' B
_{\eta(t)}\bigr)\bigr) \\
&&{}    + \bigl(a\bigl(\bar Y'_{\eta(t)}\bigr)-a\bigl(\bar\ups'_{\eta(t)}\bigr)\bigr) \bigl(\Sig' B
_t-\Sig' B _{\eta(t)}\bigr).
\end{eqnarray*}
Similar as in the proof of Proposition~\ref{prop0306-2}, we apply
Lemma \ref{le1121-1} to deduce that
%
\begin{eqnarray}\label{eq0218-1}
&&\IE [\|\bar Y-\bar Y'-(\bar\ups-\bar\ups')\|^2 ]^{1/2}\nonumber\\
&&\qquad
\le K   \IE[(\|\bar Y'\|+1)^2]^{1/2}  \IE\Bigl[\sup_{t\in
[0,1]}\bigl|L_t''-L''_{\eta(t)}-\bigl(\Sig' B _t-\Sig' B _{\eta
(t)}\bigr)\bigr|^2\Bigr]^{1/2}\\
&&{} \quad \qquad   +K  \IE[\|\bar Y'-\bar\ups'\|^2]^{1/2}   \IE\Bigl[\sup_{t\in
[0,1]} \bigl|\Sig' B _t-\Sig' B _{\eta(t)}\bigr|^2\Bigr]^{1/2}.\nonumber
\end{eqnarray}

Next, we estimate $\IE[\sup_{t\in[0,1]}|L_t''-L''_{\eta(t)}-(\Sig'
B _t-\Sig' B _{\eta(t)})|^2]$. Recall that conditional on $\IIJ$,
each pairing of $(L''_{t+T_j}-L''_{T_j})_{t\in[0,T_{j+1}-T_j]}$ and
$(B _{t+T_j}-B _{T_j})_{t\in[0,T_{j+1}-T_j]}$ is coupled according to
Theorem~\ref{th0217-1}, and individual pairs are independent of each other.

Let us first assume that the times in $\IIJ$ are deterministic with
mesh smaller or equal to $\eps'$. We denote by $n$ the number of
entries of $\IIJ$ which fall into $[0,1]$, and we denote,
for $j=1,\ldots,n$, $\Delta_j=\sup_{t\in[T_{j-1},T_{j}]}
|L_t''-L''_{T_{j-1}}-(\Sig' B _t-\Sig' B _{T_{j-1}})|$. By (\ref
{eq0325-1}) and the Markov inequality, one has, for $u\ge0$,
\[
\IP\Bigl(\sup_{j=1,\ldots,n} \Delta_j\ge u\Bigr) \le\sum_{j=1}^n \IP(\Delta
_j\ge u) \le n \exp \biggl\{ c_2\log \biggl(\frac{F(h) \eps
'}{h^2}\vee e \biggr) -\frac{c_1}{\sqrt\vth h} u \biggr\}.
\]
Let now $\alpha=\frac{c_1}{\sqrt\vth h}$, $\beta=\frac
{F(h)}{h^2}$  and $u_0=\sfrac1\alpha(\log n+c_2 \log(\beta\eps'
\vee e))$. Then for $u\ge0$
\[
\IP\Bigl(\sup_{j=1,\ldots,n} \Delta_j\ge u\Bigr)\le e^{-\alpha(u-u_0)}
\]
so that
\begin{eqnarray*}
\IE \Bigl[\sup_{j=1,\ldots,n} \Delta_j^2 \Bigr]&=&2 \int_0^\infty u
\IP\Bigl( \sup_{j=1,\ldots,n} \Delta_j\ge u\Bigr)\,\dd u\\
&\le& u_0^2 +2\int_{u_0}^\infty e^{-\alpha(u-u_0)} \,\dd u= u_0^2
+2\frac1\alpha u_0 +2\frac1{\alpha^2}\le \biggl(u_0+\frac2\alpha
 \biggr)^2.
\end{eqnarray*}
Note that the upper bound depends only on the number of entries in
$\IIJ\cap[0,1]$, and, since $\#(\IIJ\cap[0,1])$ is uniformly
bounded by $\frac2 {\eps'}+1$, we thus get in the general random
setting that
\begin{eqnarray*}
&&\IE \Bigl[\sup_{t\in[0,1]}\bigl|L_t''-L''_{\eta(t)}-\bigl(\Sig' B _t-\Sig'
B _{\eta(t)}\bigr)\bigr|^2 \Bigr]^{1/2}\\
&&\qquad\le\frac{\sqrt{\vth} h}{c_1} \biggl[
\log \biggl(1+\frac2{\eps'} \biggr) + c_2 \log\biggl(\frac{F(h)\eps
'}{h^2}\vee e\biggr)+2 \biggr].
\end{eqnarray*}
Together with Lemma \ref{le0305-1}, this gives the appropriate upper
bound for the first summand in (\ref{eq0218-1}).

By the argument preceding (\ref{eq0512-1}), one has
\[
\IE \Bigl[\sup_{t\in[0,1]} \bigl|\Sig' B _t-\Sig' B _{\eta(t)}\bigr|^2
\Bigr]^{1/2}\le\kappa_1 |\Sig'|  \sqrt{\eps'\log\frac e{\eps'}}=
\kappa_1   \sqrt{F(h) \eps'\log\frac e{\eps'}},
\]
where $\kappa_1$ is a constant that depends only on $d_X$. This
estimate is used for the second summand in~(\ref{eq0218-1}) and
putting everything together yields the statement.
\end{pf}

\section{Proof of the main results}
\label{sec:proofs}

\subsection*{Proof of Theorem~\protect\ref{th:le}}
We consider a multilevel Monte Carlo algorithm $\widehat S\in\cA$
partially specified by $\eps_k:=2^{-k}$ and $h_k:=g^{-1}(2^{k})$ for
$k\in\IZ_+$. The maximal index $m\in\IN$ and the number of
iterations $n_1,\ldots,n_{m}\in\IN$ are fixed explicitly below in
such a way that $h_m\le\mathfrak{h}$ and $m\ge2$.
Recall that
\[
\operatorname{mse}(\widehat S)\le\cW\bigl(Y,\ups^{({m})}\bigr)^2 + \sum
_{k=2}^{m} \frac1{n_k} \IE\bigl[\bigl\| \ups^{({k})}-\ups^{({k-1})}\bigr\|
^2\bigr] +\frac1{n_1} \IE\bigl[\bigl\| \ups^{({1})}-y_0\bigr\|^2\bigr];
\]
see (\ref{eq0311-1}). We control the Wasserstein metric via
Corollary~\ref{cor0311-1}. Moreover, we deduce from \cite{DerHei10}, Theorem~2, that there exists a constant $\kappa_0$ that depends only
on $K$ and $d_X$ such that, for $k=2,\ldots,m$,
\[
\IE\bigl[\bigl\| \ups^{({k})}-\ups^{({k-1})}\bigr\|^2\bigr]\le\kappa_0  \bigl(\eps
_{k-1} \log(e/\eps_{k-1})+F(h_{k-1})\bigr)
\]
 and
 \[
 \IE\bigl[\bigl\| \ups
^{({1})}-y_0\bigr\|^2\bigr]\le\kappa_0  \bigl(\eps_0 \log(e/\eps_0)+F(h_0)\bigr).
\]

Consequently, one has
%
\begin{equation}\label{eq0313-1}\qquad
\operatorname{mse}(\widehat S) \le\kappa_1  \biggl[ \biggl(h_m^2 \frac1{\sqrt
{\eps_m}} +\eps_m\biggr) \log\frac e{\eps_m}+ \sum_{k=0}^{m-1} \frac
1{n_{k+1}}  \biggl[ F(h_k)+ \eps_{k} \log\frac e{\eps_{k}}
\biggr] \biggr]
\end{equation}
in the general case, and
%
\begin{eqnarray}\label{eq0313-2}
\operatorname{mse}(\widehat S) &\le&\kappa_2 \Biggl [ h_m^2 \frac1{\sqrt
{\eps_m}} \log\frac e{\eps_m} +|b-F_0(h)|^2\eps_m^2\nonumber
\\[-8pt]
\\[-8pt]&&\hspace*{19pt}{} + \sum
_{k=0}^{m-1} \frac1{n_{k+1}}  \biggl[ F(h_k)+ \eps_k \log\frac e{\eps
_k}  \biggr] \Biggr]\nonumber
\end{eqnarray}
in the case where $\Sig=0$.
Note that $F(h_k) \le h_k^2 g(h_k)=g^{-1}(2^{k})^2 2^{k}$. With Lemma~\ref{le0313-1}, we conclude that $h_k=g^{-1}(2^k)\succsim(\gamma
/2)^k$ so that $\eps_k\log\frac e{\eps_k}=2^{-k} \log(e2^k)\precsim
g^{-1}(2^k)^2 2^k$. Hence, we can bound $F(h_k)+ \eps_k \log\frac
e{\eps_k}$ from above by a multiple of $h_k^2 g(h_k)$ in (\ref
{eq0313-1}) and (\ref{eq0313-2}).

By Lemma~\ref{le0313-1}, we have $|F_0(h_m)|\precsim h_m/\eps_m$ as
$m\to\infty$. Moreover, in the case with general $\Sigma$ and
$g^{-1}(x)\succsim x^{-3/4}$, we have $h_m^2 \frac1{\sqrt{\eps_m}}
\succsim\eps_m$. Hence, in case (I), there exists a constant $\kappa
_3$ such that
%
\begin{equation}\label{eq0313-4}
\operatorname{mse}(\widehat S) \le\kappa_3  \Biggl[ h_m^2 \frac1{\sqrt
{\eps_m}} \log\frac e{\eps_m} + \sum_{k=0}^{m-1} \frac1{n_{k+1}}
h_k^2  g(h_k)  \Biggr].
\end{equation}
Conversely, in case (II), i.e. $g^{-1}(x)\precsim x^{-3/4}$, the term
$h_m^2 \frac1{\sqrt{\eps_m}}$ is negligible in~(\ref{eq0313-1}),
and we get
%
\begin{equation}\label{eq0313-3}
\operatorname{mse}(\widehat S) \le\kappa_4 \Biggl [ \eps_m \log\frac
e{\eps_m}+ \sum_{k=0}^{m-1} \frac1{n_{k+1}} h_k^2  g(h_k) \Biggr]
\end{equation}
for an appropriate constant $\kappa_4$.

Now, we specify $n_1,\ldots,n_{m}$ in dependence on a positive
parameter $Z$ with $Z\ge1/g^{-1}(2^m)$. We set
$n_{k+1}=n_{k+1}(Z)=\lfl Z g^{-1}(2^k)\rfl\ge\frac12 Z g^{-1}(2^k)$
for $k=0,\ldots,m-1$ and conclude that, by (\ref{eq1126-4}),
%
\begin{eqnarray}\label{eq0313-5}
\sum_{k=0}^{m-1} \frac1{n_{k+1}} h_k^2  g(h_k)
&= &\sum_{k=0}^{m-1}
\frac1{n_{k+1}} 2^k g^{-1}(2^k)^2 \le\kappa_5 \frac1Z\sum
_{k=0}^{m-1} 2^k g^{-1}(2^m)  \biggl(\frac2\gamma \biggr)^{m-k}\nonumber \\
&= &\kappa_5 \frac1Z 2^m g^{-1}(2^m) \sum_{k=0}^{m-1} \gamma
^{-(m-k)}\\
&\le&\kappa_5 \frac1{1-\gamma^{-1}} \frac1{Z} 2^m g^{-1}(2^m).\nonumber
\end{eqnarray}
Similarly, we get with (\ref{eq0514-2})
%
\begin{equation}\label{eq0312-2}
\operatorname{cost}(\widehat S)\le3\sum_{k=0}^{m-1} 2^{k+1} n_k\le\kappa
_6 Z 2^m g^{-1}(2^m).
\end{equation}

We proceed with case (I). By (\ref{eq0313-4}) and (\ref{eq0313-5}),
%
\begin{equation}\label{eq0312-1}
\operatorname{mse}(\widehat S) \le\kappa_7 \biggl[ g^{-1}(2^m)^2 2^{m/2} m
+ \frac1Z 2^m g^{-1}(2^m) \biggr]
\end{equation}
so that, for $Z:=2^{m/2}/({m g^{-1}(2^m)})$,
\[
\operatorname{mse}(\widehat S) \le2 \kappa_7 g^{-1}(2^m)^2 2^{m/2} m
\]
and, by (\ref{eq0312-2}),
\[
\operatorname{cost}(\widehat S)\le\kappa_6 \frac{2^{ (3/2) m}}{m}.
\]

For a positive parameter $\tau$, we choose $m=m(\tau)\in\IN$ as the
maximal integer with $\kappa_6 2^{ (3/2) m}/m\le\tau$. Here, we
suppose that $\tau$ is sufficiently large to ensure the existence of
such a $m$ and the property $h_m\le\mathfrak h$.
Then
$\operatorname{cost}(\widehat S)\le\tau$. Since $2^m\approx(\tau\log
\tau)^{2/3}$, we conclude that
\[
\operatorname{mse}(\widehat S) \precsim g^{-1} ((\tau\log\tau
)^{2/3} )^2 \tau^{1/3} (\log\tau)^{4/3}.
\]

It remains to consider case (II). Here, (\ref{eq0313-3}) and (\ref
{eq0313-5}) yield
\[
\operatorname{mse}(\widehat S) \le\kappa_8\biggl [ 2^{-m} m + \frac1Z 2^m
g^{-1}(2^m) \biggr]
\]
so that, for $Z:=\frac1m 2^{2m} g^{-1}(2^m)$,
\[
\operatorname{mse}(\widehat S) \le2 \kappa_8 2^{-m} m
\]
and, by (\ref{eq0312-2}),
\[
\operatorname{cost}(\widehat S)\le\kappa_6 \frac1m 2^{3m} g^{-1}(2^m)^2.
\]
Next, let $l\in\IN$ such that $ 2 \kappa_6 2^{-l}\gamma^{-2l}\le
1$. Again we let $\tau$ be a positive parameter which is assumed to be
sufficiently large so that we can pick $m=m(Z)$ as the maximal natural
number larger than $l$ and satisfying $2^{m+l}\le g^{*}(\tau)$. Then,
by~(\ref{eq1126-3}),
\[
\operatorname{cost} (\widehat S)\le\kappa_6 \frac1m 2^{3m}
g^{-1}(2^m)^2\le2 \kappa_6 2^{-3l}  \biggl(\frac2\gamma \biggr)^{2l}
\frac1{m+l} 2^{3(m+l)} g^{-1}(2^{m+l})^2\le\tau.
\]
Conversely, since $2^{-m} \le2^{l+1} g^*(\tau)$,
\[
\operatorname{mse}(\widehat S) \le2 \kappa_8 2^{l+1} g^*(\tau)^{-1} \log
_2 g^*(\tau).
\]
Moreover, $g^{-1}(x)\succsim x^{-1}$ so that $
x^3 g^{-1}(x)^2/\log x \succsim x /\log x$, as $x\to\infty$. This
implies that
$\log g^*(\tau)\precsim\log\tau$.

\subsection*{Proof of Corollary~\protect\ref{cor1}}
We fix $\beta'\in(\beta,2]$ or $\beta'=2$ in the case where $\beta
=2$, and note that, by definition of $\beta$,
\[
\kappa_1:=\int_{B(0,1)}|x|^{\beta'} \nu(\dd x)
\]
is finite. We consider $\bar g\dvtx (0,\infty)\to(0,\infty), h\mapsto
\int\frac{|x|^2}{h^2}\wedge1 \nu(\dd x)$. For $h\in(0,1]$, one has
\begin{eqnarray*}
\bar g(h)&=& \int_{B(0,1)} \frac{|x|^2}{h^2}\wedge1 \nu(\dd x)+\int
_{B(0,1)^c} \frac{|x|^2}{h^2}\wedge1 \nu(\dd x)\\
&\le&
\int_{B(0,1)}\frac{|x|^{\beta'}} {h^{\beta'}} \nu(\dd
x)+\int_{B(0,1)^c} 1  \nu(\dd x) \le\kappa_2 h^{-\beta'},
\end{eqnarray*}
where $\kappa_2=\kappa_1+\nu(B(0,1)^c)$.
Hence, we find a decreasing and invertible function $g:(0,\infty)\to
(0,\infty)$ that dominates $\bar g$ and satisfies $g(h)=\kappa_2
h^{-\beta'}$ for $h\in(0,1]$. Then for $\gamma=2^{1-1/\beta'}$, one has
$g(\frac\gamma2 h)=2 g(h)$ for $h\in(0,1]$ and we are in the
position to apply Theorem~\ref{th:le}: In the first case, we get
\[
\operatorname{err}(\tau)\precsim\tau^{- ({4-\beta'})/({6\beta'})} (\log
\tau)^{ (2/3) (1- 1/{\beta'})}.
\]
In the second case, we assume that $\beta'\le\frac43$ and obtain
$g^*(\tau)\approx(\tau\log\tau)^{- {\beta'}/({3\beta'-2})}$ so that
\[
\operatorname{err}(\tau)\precsim\tau^{- {\beta'}/({6\beta'-4})} (\log
\tau)^{ ({\beta'-1})/({3\beta'-2})}.
\]
These estimates yield immediately the statement of the corollary.

\begin{appendix}
\section*{Appendix}

\begin{lemma}\label{le0713-1}
Let $(A_t)$ be a previsible process with state space $\IR^{d_Y\times
d_X}$, let $(L_t)$ be a square integrable $\IR^{d_X}$-valued L\'evy
martingale and denote by $\langle L\rangle$ the process given via
\[
\langle L\rangle_t= \sum_{j=1}^{d_X}\bigl\langle L^{({j})}\bigr\rangle_t,
\]
where $\langle L^{({j})}\rangle$ denotes the predictable compensator
of the classical bracket process for the $j$th coordinate of $L$. One
has, for any stopping time $\tau$ with finite expectation $\IE\int
_0^\tau|A_s|^2\,\dd\langle L\rangle_s$, that $(\int_0^{t\wedge\tau
} A_s\,\dd L_s)_{t\ge0}$ is a uniformly square integrable martingale
which satisfies
\[
\IE \biggl|\int_0^\tau A_s\,\dd L_s  \biggr|^2\le\IE\int_0^\tau
|A_s|^2\,\dd\langle L\rangle_s.
\]
\end{lemma}

The statement of the lemma follows from the It\^o isometry for L\'evy
driven stochastic differential equations. See, for instance, \cite
{DerHei10}, Lemma~3, for a proof.

\begin{lemma}\label{le0305-1}
The processes $\bar Y'$ and $\ups$ introduced in Section~\ref{sec31} satisfy
\[
\IE \Bigl[\sup_{s\in[0,1]} |\bar Y'_s-y_0| \Bigr] \le\kappa
\quad \mbox{and}\quad     \IE \Bigl[\sup_{s\in[0,1]} |\bar\ups
_s-y_0| \Bigr] \le\kappa,
\]
where $\kappa$ is a constant that depends only on $K$.
\end{lemma}

\begin{pf}
The result is proven via a standard Gronwall inequality type argument
that is similar to the proofs of the above propositions. It is
therefore omitted.
\end{pf}

\begin{lemma}\label{le0313-1}
Let $\bar h>0$, $\gamma\in(1,2)$  and $g\dvtx (0,\infty)\to(0,\infty)$
be an invertible and decreasing function such that, for $h\in(0,\bar h]$,
\[
g\biggl(\frac\gamma2 h\biggr)\ge2 g(h).
\]
Then
%
\begin{equation}\label{eq1126-3}
\frac\gamma2 g^{-1}(u)\le g^{-1}(2u)
\end{equation}
for all $u\ge g(\bar h)$.
Moreover, there exists a finite constant $\kappa_1$ depending only on
$g$ such that for all $k,l\in\IZ_+$ with $k\le l$ one has
%
\begin{equation}\label{eq1126-4}
g^{-1}(2^k)\le\kappa_1 \biggl (\frac{2}\gamma \biggr)^{l-k} g^{-1}(2^l).
\end{equation}

If $\nu(B(0,h)^c)\le g(h)$ for all $h>0$, and $\nu$ has a second
moment, then
\[
\int_{B(0,h)^c} |x| \nu(\dd x) \le\kappa_2\bigl(h g(h)+1\bigr),
\]
where $\kappa_2$ is a constant that depends only on $g$ and $\int
|x|^2 \nu(\dd x)$.
\end{lemma}

\begin{pf}
First, note that property (\ref{eq1126-1}) is equivalent to
\[
\frac\gamma2 g^{-1}(u)\le g^{-1}(2u)
\]
for all sufficiently large $u>0$. This implies that there exists a
finite constant $\kappa_1$ depending only on $g$ such that for all
$k,l\in\IZ_+$ with $k\le l$ one has
\[
g^{-1}(2^k)\le\kappa_1  \biggl(\frac{2}\gamma \biggr)^{l-k} g^{-1}(2^l).
\]
For general, $h>0$ one has
\[
\int_{B(0,h)^c} |x| \nu(\dd x) \le\int_{B(0,h)^c \cap B(0,\bar h)}
|x| \nu(\dd x) +\frac1{\bar h} \int|x|^2 \nu(\dd x).
\]
Moreover,\vspace*{-1pt}
\begin{eqnarray*}
\int_{B(0,h)^c \cap B(0,\bar h)} |x| \nu(\dd x)
&\le&\sum
_{n=0}^\infty\nu \biggl(B\biggl(0,h\biggl(\sfrac2\gamma\biggr)^n \biggr)^c\cap B(0,\bar
h)\biggr) h \biggl(\sfrac2\gamma \biggr)^{n+1}\\
&\le&\sum_{n=0}^\infty\ind_{\{h( 2/\gamma)^n\le\bar h\}}
\underbrace{g \biggl(h\biggl(\sfrac2\gamma\biggr)^n \biggr)}_{\le2^{-n} g(h)} h
 \biggl(\sfrac2\gamma \biggr)^{n+1}\\[-1pt]
&\le&2h g(h)\sum_{n=0}^\infty\gamma^{-(n+1)}.
\end{eqnarray*}
\upqed
\end{pf}

\begin{lemma}\label{le1121-1} Let $n\in\IN$ and $(\cG
_j)_{j=0,1,\ldots,n}$ denote a filtration.
Moreover, let, for $j=0,\ldots,n-1$, $U_j$ and $V_j$ denote nonnegative
random variables such that
$U_j$ is $\cG_{j}$-measurable, and $V_j$ is $\cG_{j+1}$-measurable
and independent of $\cG_{j}$. Then one has\vspace*{-1pt}
\[
\IE \Bigl[\max_{j=0,\ldots,n-1} U_jV_j  \Bigr]\le\IE \Bigl[\max
_{j=0,\ldots,n-1}U_j  \Bigr]\cdot\IE \Bigl[\max_{j=0,\ldots,n-1}
V_j \Bigr].
\]
\end{lemma}
\begin{pf}
See \cite{DerHei10}.
\end{pf}

\end{appendix}

\printaddresses

\end{document}